\documentclass{amsart}

\title{The contingent HOD dichotomy}
\usepackage[pdfauthor={Joel David Hamkins, Bokai Yao},
    pdftitle={The maximality principle and V=HOD},
    hidelinks]{hyperref}

\author{Joel David Hamkins}
\address[Joel David Hamkins]
{O'Hara Professor of Logic, University of Notre Dame, 100 Malloy Hall, Notre Dame, IN 46556 USA, \&\ 
Guest Chair Professor, Peking University, Beijing, China}
\email{jdhamkins@nd.edu}
\urladdr{http://jdh.hamkins.org}
\thanks{This paper grew out of a question of the second author that arose during the first author's lectures on the modal logic of forcing at Peking University in June 2026. The first author is grateful for the support of a Guest Chair Professorship provided by Peking University. A preliminary form of the contingent \HOD\ dichotomy was proved by the first author on his blog~\cite{Hamkins2015:HOD-of-a-set}. The authors are grateful to Gunter Fuchs for observations concerning the blurry \HOD\ dichotomy in connection with the contingent \HOD\ dichotomy. The forcing modalities were introduced by the first author~\cite{Hamkins2003:MaximalityPrinciple} and further investigated in collaborative work on the modal logic of forcing~\cite{HamkinsLoewe2008:TheModalLogicOfForcing, HamkinsLeibmanLoewe2015:StructuralConnectionsForcingClassAndItsModalLogic}. The maximality principles were introduced and investigated by Stavi and \Vaananen~\cite{StaviVaananen2001:ReflectionPrinciples} 
and independently, with the forcing modalities, by the first author~\cite{Hamkins2003:MaximalityPrinciple}. The authors have a followup project underway with Yuanshan (Tom) Li on the \ZF\ analogue of the contingent \HOD\ dichotomy.}

\author{Bokai Yao}
\address[Bokai Yao]
{Department of Philosophy and Religious Studies, Peking University, China}
\email{bkyao@pku.edu.cn}
\urladdr{bokaiyao.com}

\usepackage{latexsym,amsfonts,amsmath,amssymb,mathrsfs}
\usepackage{tikz}
\usetikzlibrary{arrows,arrows.meta,petri,topaths,positioning,shapes,shapes.misc,patterns,calc,decorations.pathreplacing,hobby}
\usepackage[shortlabels]{enumitem} 
\usepackage{harpoon}

\usepackage{bbm}
\usepackage[backend=biber,style=alphabetic,maxbibnames=15,maxcitenames=6,dateabbrev=false]{biblatex}
\renewcommand{\UrlFont}{} 
\renewbibmacro{in:}{\ifentrytype{article}{}{\printtext{\bibstring{in}\intitlepunct}}} 
\DeclareFieldFormat{url}{\UrlFont\url{#1}} 
\DeclareFieldFormat{urldate}{
  (version \thefield{urlday}\addspace%
  \mkbibmonth{\thefield{urlmonth}}\addspace%
  \thefield{urlyear}\isdot)}
\DeclareFieldFormat{eprint:arxiv}{
  \ifhyperref
    {\href{http://arxiv.org/abs/#1}{%
        arXiv\addcolon\nolinkurl{#1}}\iffieldundef{eprintclass}{}{\UrlFont{\mkbibbrackets{\thefield{eprintclass}}}}}
    {arXiv\addcolon\nolinkurl{#1}\iffieldundef{eprintclass}{}{\UrlFont{\mkbibbrackets{\thefield{eprintclass}}}}}}
\newtheorem{theorem}{Theorem}
\newtheorem*{theorem*}{Theorem}

\newtheorem*{maintheorem*}{Main Theorem}
\newtheorem*{maintheorems*}{Main Theorems}
\newtheorem{corollary}[theorem]{Corollary}
\newtheorem*{corollary*}{Corollary}
\newtheorem*{corollaries*}{Corollaries}

\newtheorem{observation}[theorem]{Observation}

\theoremstyle{definition}

\newtheorem*{definition*}{Definition}

\newtheorem{question}[theorem]{Question}
\newtheorem*{question*}{Question}

\newtheorem*{questions*}{Questions}
\newtheorem*{mainquestion*}{Main Question} 
\newtheorem*{openquestion*}{Open Question} 
\theoremstyle{remark}

\theoremstyle{plain}
\newcommand{\one}{\mathbbm{1}} 
\newcommand{\set}[1]{\{\,{#1}\,\}}
\newcommand{\of}{\subseteq}
\newcommand{\ofneq}{\subsetneq}
\newcommand{\MP}{{\rm MP}}

\newcommand\blurryHOD{\text{\kern1pt-}\HOD}
\newcommand{\UnderTilde}[1]{{\setbox1=\hbox{$#1$}\baselineskip=0pt\vtop{\hbox{$#1$}\hbox to\wd1{\hfil$\sim$\hfil}}}{}}
\newcommand{\HOD}{{\rm HOD}}
\newcommand{\HOA}{{\rm HOA}}
\newcommand{\OD}{{\rm OD}}
\newcommand{\AC}{{\rm AC}}
\newcommand{\Con}{\mathop{{\rm Con}}}
\newcommand{\satisfies}{\models}
\newcommand{\forces}{\Vdash}

\DeclareMathOperator{\possible}{\text{\tikz[scale=.6ex/1cm,baseline=-.6ex,rotate=45,line width=.1ex]{\draw (-1,-1) rectangle (1,1);}}}
\DeclareMathOperator{\necessary}{\text{\tikz[scale=.6ex/1cm,baseline=-.6ex,line width=.1ex]{\draw (-1,-1) rectangle (1,1);}}}
\newcommand{\CCA}{{\rm CCA}}
\newcommand{\CH}{{\rm CH}}
\newcommand{\GCH}{{\rm GCH}}
\newcommand{\SCH}{{\rm SCH}}
\newcommand{\ZFC}{{\rm ZFC}}
\newcommand{\ZF}{{\rm ZF}}
\newcommand{\KM}{{\rm KM}}

\newcommand{\AD}{{\rm AD}}
\newcommand{\Tr}{\mathop{\rm Tr}\nolimits}
\newcommand{\elesub}{\prec}

\newcommand{\Union}{\bigcup}
\renewcommand{\implies}{\mathrel{\rightarrow}}
\newcommand{\Implies}{\mathrel{\Longrightarrow}}

\renewcommand{\iff}{\mathrel{\leftrightarrow}}

\newcommand{\M}{{\mathbb M}}
\renewcommand{\P}{{\mathbb P}}
\newcommand{\Q}{{\mathbb Q}}

\newcommand{\R}{{\mathbb R}}
\def\<#1>{\left\langle#1\right\rangle}
\newcommand\Vaananen{Väänänen}
\newcommand{\Add}{\mathop{\rm Add}}
\newcommand{\GA}{{\rm GA}}
\newcommand{\Ord}{\mathord{{\rm Ord}}}

\newcommand{\Qdot}{{\dot\Q}}

\newcommand{\Levy}{L\'{e}vy}
\newcommand\Vopenka{Vop\v{e}nka}
\newcommand{\Bukovsky}{Bukovsk\`y}
\newbox\gnBoxA
\newbox\gnBoxB
\newdimen\gnCornerHgt
\setbox\gnBoxA=\hbox{\tiny$\ulcorner$}
\global\gnCornerHgt=\ht\gnBoxA
\newdimen\gnArgHgt
\def\gcode #1{%
\setbox\gnBoxA=\hbox{$#1$}%
\setbox\gnBoxB=\hbox{$\bar #1$}%
\gnArgHgt=\ht\gnBoxB%
\ifnum     \gnArgHgt<\gnCornerHgt \gnArgHgt=0pt%
\else \advance \gnArgHgt by -\gnCornerHgt%
\fi \raise\gnArgHgt\hbox{\tiny$\ulcorner$} \box\gnBoxA %
\raise\gnArgHgt\hbox{\tiny$\urcorner$}}
\newcommand{\restrict}{\upharpoonright} 
\newcommand*{\vect}[1]{\overrightharp{\ensuremath{#1}}}
\renewcommand{\emptyset}{\varnothing}
\newcommand{\Mantle}{\mathbb{M}}
\newcommand{\gHOD}{\ensuremath{\mathord{{\rm g}\HOD}}} 
\newcommand{\intersect}{\cap}
\newcommand{\Intersect}{\bigcap}

\begin{document}

\begin{abstract}
The axiom $V=\HOD$, referring to the class of hereditarily ordinal-definable sets, is contingently contingent by forcing---it holds contingently in the generic multiverse of some models of set theory but not others. Thus we introduce and investigate the \emph{contingent \HOD\ dichotomy}, according to which either (1) the set-theoretic universe $V$ is close to \HOD\ by forcing---\HOD\ is a ground, $V=\HOD$ is forceable, the universe is ordinal-definable from a set parameter $V=\HOD_b$, the blurry \HOD\ hierarchy collapses to $\kappa\blurryHOD$ for some cardinal $\kappa$, and furthermore these features are invariant by forcing and hold throughout the generic multiverse; or (2) the universe $V$ is far from \HOD\ by forcing in all those respects, so $V\neq\HOD$, but also \HOD\ is not a ground, $V=\HOD$ is not forceable, the set-theoretic universe is not ordinal-definable from a set, the universe is not $\kappa\blurryHOD$ for any $\kappa$, and these features hold invariantly throughout the generic multiverse. We investigate the nature of the dichotomy in connection with large cardinals, with set-theoretic geology, with the iterated $\HOD^\alpha$ classes, and with the maximality principles $\MP$ and $\MP(\R)$, showing generally that both prongs of the dichotomy occur with diverse set-theoretic phenomena.
\end{abstract}

\maketitle

\section{Introduction}

The \emph{contingent \HOD\ dichotomy} expresses the two distinct ways in which the set-theoretic universe $V$ in \ZFC\ can relate to the class of hereditarily ordinal definable sets by forcing. Namely, either (1) the universe $V$ is close to \HOD\ by forcing in that \HOD\ is a ground model, the axiom $V=\HOD$ is forceable, the universe is ordinal-definable from a set parameter, the blurry \HOD\ hierarchy collapses, so $V=\kappa\blurryHOD$ for some cardinal $\kappa$, and furthermore each of these things is invariant by forcing to all forcing extensions and grounds, being true throughout the generic multiverse; or (2) the universe $V$ is far from \HOD\ by forcing, in particular $V\neq\HOD$, but also \HOD\ is not a ground, the axiom $V=\HOD$ is not forceable, the universe is not ordinal-definable from a parameter, the blurry \HOD\ hierarchy does not collapse, so $V\neq\kappa\blurryHOD$ for any $\kappa$, and furthermore these things are invariant by forcing throughout the generic multiverse. The first prong of the dichotomy is thus characterized by a contingent nature of \HOD\ with respect to forcing---you can turn it on or off as desired---whereas in the second prong we find a necessary failure $\necessary(V\neq\HOD)$ throughout the generic multiverse. If \ZFC\ is consistent, then both prongs of the dichotomy occur.

In this paper, we shall introduce and establish the dichotomy and investigate how it is manifested in regard to other set-theoretic phenomena, such as large cardinals, the \GCH, set-theoretic geology, the iterated $\HOD^\alpha$ classes, and the maximality principles in the modal logic of forcing.

The contingent \HOD\ dichotomy is similar to but not identical with the \HOD\ dichotomy of W.~Hugh Woodin \cite{Woodin_Davis_Rodriguez_2012}. To be sure, both dichotomies are accurately described as a distinction between $V$ being ``close to \HOD'' or ``far from \HOD'' in robust relevant senses, so the dichotomies share a natural affinity. Whereas Woodin's \HOD\ dichotomy is aimed at a deep Jensen-covering phenomenon regarding \HOD, the contingent \HOD\ dichotomy, in contrast, is concerned with softer forcing-invariant features relating the universe to \HOD. It happens that the first prong of our contingent \HOD\ dichotomy here always lies within the first prong of Woodin's \HOD\ dichotomy, but ultimately the two dichotomies do not draw exactly the same line. The contingent \HOD\ dichotomy thus in effect refines what is after all the main prong of Woodin's dichotomy, as he conjectures that his second prong does not occur in the relevant large cardinal context. We discuss the matter further in section \ref{Section.Contingent-HOD-dichotomy-versus-HOD-dichotomy}.

We shall begin with a review of \HOD, including a discussion of two sometimes puzzling elementary matters concerning the definition and history of \HOD\ before moving on to the core contingency features. We shall prove the contingent \HOD\ dichotomy in section \ref{Section.The-contingent-HOD-dichotomy} and then move on to explore its relation with other set-theoretic phenomena in the subsequent sections.

\section{Review of $\HOD$ and a certain philosophical puzzle}

Let us quickly review the basics of $\HOD$, the class of hereditarily ordinal-definable sets. The axiom $V=\HOD$ asserts that every object in the set-theoretic universe is definable from ordinal parameters. 
That is, for every set $a$ there is a formula $\varphi$ in the language of set theory and ordinal parameters $\vect\alpha$, such that $a$ is the unique object satisfying $\varphi(a,\vect\alpha)$. To help in the analysis, set theorists define the associated class \OD\ consisting of all \emph{ordinal-definable} sets, the sets $a$ for which there is such a definition $\varphi(a,\vect\alpha)$.

The class of ordinal-definable sets $\OD$ is not necessarily transitive---a set could be ordinal definable without all its elements being ordinal definable. For example, the set of real numbers is definable (using any of several standard accounts), even if some individual real numbers are not ordinal definable, as is known to occur in some models of set theory. To find the transitive inner core of \OD, therefore, set theorists define the class of \emph{hereditarily ordinal-definable} sets, denoted $\HOD$, which are the sets that are ordinal definable along with all their hereditary members---the set and everything in its transitive closure is ordinal definable. In any model of \ZF, the class $\HOD$ is a transitive inner model of \ZFC\ containing all the ordinals. This class is therefore often useful for metamathematical results in set theory. One can use $\HOD$ instead of the constructible universe $L$, for example, to prove the relative consistency of the axiom of choice. 

There is a certain subtle philosophical puzzle at the core of the definition of ordinal definability. Namely, how is it possible in the first place to define ordinal definability in the first-order language of set theory? Have we actually defined \OD\ and \HOD\ as classes in \ZFC\ set theory? In the account above, we quantified over the formulas $\varphi$ and referred to the satisfaction of $\varphi(a,\vect\alpha)$, whereas Tarski proved that there is no first-order definable truth predicate in the language of set theory, so how can we express that idea in the language of set theory? 

Indeed, we are not generally able to express ``$x$ is definable'' in set theory. A little contemplation will show that we cannot refer unproblematically to ``the smallest undefinable ordinal,'' since this notion cannot be definable if indeed there are any undefinable ordinals. Pressing further with this idea shows that in some models of \ZFC\ the collection of definable elements is not first-order definable---the definable elements do not necessarily form a \ZFC-amenable class,~\cite[p.~145]{HamkinsLinetskyReitz2013:PointwiseDefinableModelsOfSetTheory}.\footnote{Interestingly, there are other models of \ZFC\ in which the notion ``$x$ is definable'' is expressible---in any pointwise definable model, for example, every object is definable without parameters, so in such a model the class of definable elements is trivially defined by $x=x$~\cite{HamkinsLinetskyReitz2013:PointwiseDefinableModelsOfSetTheory}.}

Why are we able to express ``$x$ is ordinal definable'' when we cannot express ``$x$ is definable'' or ``$x$ is definable from parameter $p$''? In fact we \emph{are} able to express the property of being ordinal definable, uniformly so in every model of \ZF. The collection of hereditarily ordinal-definable sets $\HOD$ is a first-order definable class in \ZF\ and $V=\HOD$ is expressible as an axiom in the language of set theory.

One standard approach is to define officially that a set $x$ is \emph{ordinal-definable} if there is some ordinal $\theta$ such that $x$ is definable in $\<V_\theta,\in>$ from  ordinal parameters. That is, there is some formula $\varphi$ and ordinal parameters $\vect\alpha$ such that $x$ is the unique object in $V_\theta$ satisfying $\varphi(x,\vect\alpha)$ in that structure. This version of the concept is expressible in the language of set theory because each $V_\theta$ is a set and so we have a truth predicate satisfaction relation for the structure $\<V_\theta,\in>$ and thus are able to refer to definability with respect to truth in $V_\theta$ with ordinal parameters. What remains, of course, is to prove the meta-mathematical fact that this official definition in set theory aligns with the meta-theoretic, external conception of being ordinal definable. Some accounts of $\HOD$ introduce the Gödel operations, which are primitive functions whose values correspond at bottom to the logical connectives $\wedge,\vee,\neg$ and quantifiers $\exists,\forall$ used to build up the formulas; compositions of these operations thus correspond essentially to formulas to be applied in a given set.

Set theorists are quick to point out that if a set $a$ is definable in $V$ from ordinal parameters, in the sense that we can explicitly provide a formula $\varphi$ and ordinal parameters $\vect\alpha$ such that $a$ is the unique object for which $\varphi(a,\vect\alpha)$, then by the reflection theorem this property will reflect to some $V_\theta$, and so $a$ will be  ordinal-definable by the official definition we just gave. So this is half of the equivalence.

Conversely, one sometimes finds on offer an argument like this: if $a$ is definable in $\<V_\theta,\in>$ with formula $\varphi$ and ordinal parameters $\vect\alpha$, then we need only add $\theta$ as an additional parameter to define $a$ in $V$. Specifically, the proposal is that we may define $a$ as ``the thing for which $\varphi(a,\vect\alpha)$ holds, relativized to $V_\theta$.'' So $a$ is definable this way in $V$, using parameters $(\theta,\vect\alpha)$. Is it right?

Unfortunately, not quite. Strictly speaking this way of arguing does not always work---in fact we shall sometimes need additional ordinal parameters beyond those appearing in $\vect\alpha$ and $\theta$, as well as a slightly different defining formula. To appreciate the difficulty, notice that the official definition references definability-in-$V_\theta$ as an internal notion, a notion expressed in the object theory rather than the metatheory. When we apply the definition in an $\omega$-nonstandard model $M$, consequently, the definition includes objects $a$ that may be defined in some $V_\theta^M$ using a \emph{nonstandard} formula $\varphi$. Worse, the tuple of ordinals $\vect\alpha$ might have a nonstandard finite length, which would mean infinitely many ordinal parameters when considered outside $M$, even if internally this is considered a finite list. We cannot expect to use such a formula $\varphi$ directly or such a parameter list $(\theta,\vect\alpha)$ to show that $a$ is definable in $M$ with ordinal parameters. We want, after all, to define the object $a$ with a standard-finite formula using a standard-finite tuple of ordinals. The issue arises just as much with Gödel operations, since there will be nonstandard compositions of Gödel operations, including some with nonstandard arity.

To see that this difficulty occurs, consider an $\omega$-nonstandard model with uncountably many nonstandard natural numbers, such as in an ultrapower of the universe by a nonprincipal ultrafilter on $\omega$. This model thinks as an internal matter that every element of its (nonstandard) $V_\omega$ is definable in $V_\omega$ without parameters, since this is easily provable in \ZF\ by induction on the finite ranks. But these objects cannot all be actually definable in this $V_\omega$ without parameters, since there are only countably many definitions available in the metatheory.

Nevertheless, the problems are easily addressed, and we can prove that in any model of set theory $M$ the official internal notion of ordinal definability agrees perfectly with the outside-the-model meta-theoretic notion of being definable in $M$ from ordinal parameters. First, we can solve the issue of nonstandard finite tuples of parameters $\vect\alpha$ by using the Cantor pairing function on ordinals to prove that we only ever need at most one ordinal parameter. Every finite tuple of ordinals is canonically coded by a single ordinal, and this works perfectly well inside $M$ to code every finite tuple of ordinals, including the nonstandard finite tuples, with a single ordinal. Thus we need at most one ordinal parameter to define any ordinal-definable set. In fact, we don't need ordinal parameters $\vect\alpha$ at all beyond the ordinal $\theta$---a set is ordinal definable if and only if it is definable in some $V_\theta$ without parameters. The reason is that we may in effect code the ordinal into $\theta$ itself. If $a$ is definable in $V_\lambda$ from ordinal $\alpha$, then let $\theta=\<\lambda,\alpha>+1$, using the pairing function to code the pair, and observe that $\lambda$ and $\alpha$ are definable in $V_\theta$ without parameters, since $\<\lambda,\alpha>$ is the largest ordinal there. So $a$ is definable in $V_\theta$ without parameters. For a slick alternative proof, observe that every finite tuple of ordinals $\vect\alpha$ is definable without parameters in some $V_\theta$, for if not, there is a least counterexample, so by reflection there is some $V_\theta$ that thinks $\vect\alpha$ is the least counterexample, but this is itself a definition of $\vect\alpha$ in $V_\theta$; so there can be no counterexample.

Second, we can resolve the difficulty occurring when an element $a$ is defined in some $V_\theta^M$ using a nonstandard formula $\varphi$ by observing that the Gödel code of the formula $\gcode{\varphi}$ is a natural number in the model, even if it is nonstandard, and therefore it is available to use as an additional ordinal parameter. Consequently, if $M$ thinks $a$ is defined in $V_\theta^M$ by the nonstandard definition $\varphi(a,\alpha)$, then we can define $a$ in $M$ externally by realizing that in $M$ the set $a$ is ``the unique object for which the formula coded by $\gcode{\varphi}$ is true of $a$ in $V_\theta$ with parameter $\alpha$.'' This is a standard-finite definition using parameters $(\theta,\alpha,\gcode{\varphi})$, which are ordinals in $M$. 

Putting this all together, we conclude that the official internal definition of ordinal definable is in fact equivalent to the outside-the-universe metatheoretic notion of ordinal definable, and we may therefore freely use this concept within the language of set theory. Ordinal definability is definable, even if definability is not.

An alternative resolution of the issue proceeds like this. In $V$ we may define a global well order $\leq_{\OD}$ of the class \OD\ of ordinal definable sets, using an analogue of the global well-order in the constructible universe; this is how one knows that \AC\ holds in $\HOD$. Namely, we define that $a\leq_{\OD}b$ if $a$ is definable in a smaller $V_\theta$ than $b$, or in the same least $V_\theta$ but with a smaller formula, or with the same least $V_\theta$ and formula, but with a smaller parameter (or equal). This is a definable set-like well order of the class \OD. Every ordinal-definable set is therefore the $\alpha$th set in this order for some ordinal $\alpha$, so we can define the set this way using the definability of the order, with $\alpha$ as a parameter. So every element of \OD\ is ordinal definable in the external meta-theoretic sense, and consequently the class \OD\ coincides with the sets that are definable from ordinal parameters. So once again, we conclude that ordinal definability is definable.

Thus we see in a variety of ways that the external meta-theoretic notion of being definable from ordinal parameters is internally expressible in the language of set theory, and so we may freely use this notion in set theory. To our way of thinking, this makes for a surprising, welcome resolution of the philosophical puzzle, for it brings a concept that might otherwise have remained a philosophical meta-theoretic notion into the object theory of set theory.

\section{A historical puzzle regarding ordinal definability}

Let us also briefly discuss a certain historical puzzle about $\HOD$.\footnote{The first author is grateful to Juliette Kennedy, from whom he first heard about this particular puzzle years ago at a talk she gave in connection with her paper~\cite{Kennedy2013:On-formalism-freeness-implementing-Godel-1946}.} The concept of ordinal definability was first suggested in 1946 by Gödel in his prescient though brief remarks before the Princeton bicentennial conference on problems in mathematics. He introduces the concept as a natural transfinite instance of the process of moving to ``the next higher system,'' deducing new propositions from the resulting stronger axioms of infinity. 
\begin{quotation}\small
Here also you have, corresponding to the transfinite hierarchy of formal systems, a transfinite hierarchy of concepts of definability
\ldots 
The simplest way of doing it is to take the ordinals themselves as primitive terms. So one is led to the concept of definability in terms of ordinals, | i.e., definability by expressions containing names of ordinal numbers and logical constants, including quantification referring to sets. This concept should, I think, be investigated. It can be proved that it has the required closure property: By introducing the notion of truth for this whole transfinite language, i.e., by going over to the next language, you will obtain no new definable sets (although you will obtain new definable properties of sets). 
\cite[p.151]{Godel1946:Remarks-before-the-Princeton-bicentennial-conference-on-problems-in-mathematics}
\end{quotation}
Gödel does not elaborate with much technical detail, but clearly he proposes to express the concept of ordinal definability in terms of a formal truth predicate. Indeed, if we augment the language of set theory with a predicate for first-order truth $\Tr(\gcode{\varphi},\vect a)$ fulfilling the Tarskian recursive definition of truth, then in the expanded language we can naturally express the concept of being ordinal definable. Namely, $a$ is ordinal definable if there is a Gödel code $f$ for a formula and a tuple of ordinal parameters $\vect\alpha$ such that $a$ is the unique object for which $\Tr(f,\<a,\vect\alpha>)$. 

Gödel mentions that ``by going over to the next language, you will obtain no new definable sets,'' and we can easily prove this using the reflection theorem for the expanded language. Namely, by reflection there is an ordinal $\theta$ such that $\<V_\theta,\in,\Tr\restrict V_\theta>$ is $\Sigma_1$-elementary in $\<V,\in,\Tr>$. It follows that $V_\theta$ is closed under the existence of witnesses for the existential assertions declared true by $\Tr$, and consequently $\Tr\restrict V_\theta$ fulfills the Tarskian recursion for truth in $\<V_\theta,\in>$. Therefore truth in $V_\theta$ agrees with $\Tr$, so we can replace any instance of $\Tr(x,y)$ for $x,y\in V_\theta$ with truth in $V_\theta$, which is expressible in set theory from parameter $\theta$ without reference to $\Tr$. Thus any set that is definable from parameters in the language with the truth predicate is also definable in the language without the truth predicate, just as Gödel claimed---we need only allow $\theta$ as an additional ordinal parameter.\footnote{Kennedy describes the situation by saying ``
if one passes to the `next language,' i.e., one obtained by adding a truth predicate for statements about ordinal definable sets, one obtains no new ordinal definable sets''~\cite[p.~381]{Kennedy2013:On-formalism-freeness-implementing-Godel-1946}, adding further (fn.~139) that 
``More exactly, for any formula $\phi(x)$ there is another formula $\phi'(x)$ which says that $\phi(x)$ is true in HOD. So no new ordinal definable sets are created by referring to truth in HOD.'' This account seems to make a different interpretation of Gödel's proposal than we have, however, since we take Gödel's truth predicate to refer to truth in $V$, not $\HOD$---sets are admitted to $\HOD$ when they are (hereditarily) ordinal definable in $V$---and the relevant question concerns sets being definable in $V$, not $\HOD$.} 

The historical puzzle we want to discuss grows out of the fact that Gödel made his remarks in 1946, which was well before the reflection theorem was proved by \Levy\ and Montague~\cite{Levy1960:Axiom-schemata-of-strong-infinity-in-axiomatic-set-theory, Montague1961:Fraenkels-addition-to-axioms-of-Zermelo}. How did Gödel know how to eliminate the use of the truth predicate? Did Gödel know the reflection theorem? Did he have an alternative proof in mind?

In his introductory remarks to Gödel's article, Charles Parsons says:
\begin{quote}\small
Gödel observes that one obtains no new ordinal-definable sets by extending the language by a truth predicate, presumably for formulas containing names of arbitrary ordinals. This implies immediately that ``$x$ is an ordinal-definable set'' is expressible in the language of ZF, an assertion that has usually been proved by means of Levy's reflection principle. An extension of that argument proves Gödel's assertion about the ordinal-definable sets (OD). It is reasonable to conjecture that Gödel knew the reflection principle at the time. Otherwise he would have had to have a quite different argument in mind; moreover, an argument like the proof of the principle figures in the proof in Gödel 1940 that V = L (where V is the universe of all sets and L is the class of all constructible sets) implies the generalized continuum hypothesis.~\cite[p.~147]{Parsons1990:Introductory-note-to-1946}
\end{quote}
Gödel does not need the full reflection theorem, since as we described in our account of the argument mere $\Sigma_1$-reflection in the augmented language suffices. And surely $\Sigma_1$ reflection is easier both to imagine and establish, since one needs only to close under simple witnesses rather than undertake the more elaborate iterative argument used for full reflection---perhaps Gödel saw this much without mentioning it?

In our view the evidence that Gödel had the reflection theorem in 1946 is thin. His arguments are  brief---the main technical content appears already in the quotation we provided earlier; he does not give any details involving reflection, but simply asserts baldly that ``it can be proved that'' the expanded language offers no new definable sets; and he gives no hint that we should replace truth with true-in-a-suitable-$V_\theta$, which is the key idea at the center of the reflection argument.

We should like to emphasize that Gödel does not suggest in his 1946 remarks that the concept of ordinal definability itself is definable in the language of set theory. He simply does not offer \OD\ as a definable class in \ZF, as we now understand it to be. This observation tends to deflate some forms of the historical puzzle. If one asks, ``How was it possible for Gödel to prove, without the reflection theorem, that ordinal definability is definable in the language of set theory?'' then the answer seems to be simply that he did not prove this nor even suggest that it was true.

Gödel's proposal, after all, is to define ordinal definability in ``the next language'' with the truth predicate. Yes, he mentions that every definable set in the expanded language is definable in the language of set theory, a fact now usually proved via reflection, but this falls short of the claim that the concept of ordinal definability itself can be defined without reference to truth, a fact now also commonly proved via reflection. Indeed, he specifically warns us that the expanded language is able to define new classes, ``new definable properties of sets.'' Perhaps the clearest instance of this is the truth predicate $\Tr$ itself, which is obviously definable in the language with the truth predicate, while Tarski's theorem on the nondefinability of truth shows that $\Tr$ is not definable in the first-order language of set theory.

Curiously, Parsons (in the remarks quoted above) cites the definability of ordinal definability as an ``immediate consequence'' of the fact that there are no new definable sets. He seems to suggest that there is a simple direct reflection-free argument from no-new-definable-sets to the definability of ordinal definability, but unfortunately he offers no details and we are not aware of any such argument. What did he have in mind?

Let us mention a few further issues concerning Gödel's approach to ordinal definability. First, Gödel's truth-predicate method will not provide a definition of \OD\ in \ZF, even with the reflection theorem, since his proposal takes place in the strictly stronger theory $\ZFC(\Tr)+$``$\Tr$ is a truth predicate,'' where by $\ZFC(\Tr)$ we mean that the class $\Tr$ is \ZFC-amenable---the theory includes the collection and separation axiom schemes of \ZFC\ in the expanded language with the predicate $\Tr$. This theory altogether is strictly stronger than \ZFC\ in consistency strength and should be viewed as a weak large cardinal theory.

To see this large-cardinal nature, observe that the existence of a \ZFC-amenable truth predicate fulfilling the Tarski recursion is equivalent to realizing the universe $V$ as the union of a continuous elementary chain along a closed unbounded class of cardinals:
 $$V_{\theta_0}\elesub V_{\theta_1}\elesub V_{\theta_2}\elesub\cdots\elesub V_{\theta_\alpha}\elesub\cdots \elesub V.$$
These cardinals $\theta_\alpha$ are simply the $\Sigma_1(\Tr)$-correct cardinals, which are all worldly cardinals and indeed otherworldly to a high degree (see~\cite{Hamkins2020:Otherworldly-cardinals}), thus revealing truth-augmented set theory as a large cardinal theory. In short, the existence of a \ZFC-amenable truth predicate for first-order truth implies the existence of a proper class club of hyper-otherworldly cardinals. A similar idea is used to prove that Kelley-Morse set theory \KM\ implies $\Con(\ZFC)$ and more---\KM\ implies that there are truth predicates relative to any particular class; see~\cite{Hamkins.blog2014:KM-implies-ConZFC-and-much-more}.

Gödel's proposal that we go ``over to the next language'' is therefore not entirely innocent---one cannot just freely add a Tarskian truth predicate to a given model of set theory. If consistent, there are models of \ZFC\ that have no \ZFC-amenable truth predicate, that is, a predicate fulfilling $\ZFC(\Tr)$, which is used to prove generalized reflection. Consider a model without otherworldly cardinals. Or consider a pointwise definable model of \ZFC, one in which every element is definable without parameters. To construct such a model, take the definable elements of any model of $\ZFC+V=\HOD$; since the global well order is definable, this is closed under Skolem witnesses and is therefore an elementary submodel, hence pointwise definable since the same definitions still work. Prominent instances would include the Shepherdson-Cohen model, which is the smallest transitive model of \ZFC, but also every countable model of \ZFC\ admits a class forcing extension that is pointwise definable~\cite{HamkinsLinetskyReitz2013:PointwiseDefinableModelsOfSetTheory}, and every countable model of \ZFC\ has a pointwise definable top extension~\cite{Hamkins2024:Every-countable-model-of-arithmetic-or-set-theory-has-a-pointwise-definable-end-extension}. A pointwise definable model can have no \ZFC-amenable truth predicate, since the predicate would reveal the pointwise definable nature, causing the model to see itself as countable. In light of these models, we see that Gödel's approach to ordinal definability via truth simply does not apply in all models of \ZFC. 

Gödel had speculated that the class of ordinal definable sets would be a model of the axioms of set theory. 
\begin{quote}
In conclusion I would like to say that, irrespective of whether the concept of definability suggested in this lecture corresponds to certain intuitive notions, it has some intrinsic mathematical interest; in particular, there are two questions arising in connection with it: (1) Whether the sets definable in this sense satisfy the axioms of set theory. I think this question is to be answered in the affirmative, and so will lead to another, and probably simpler, proof for the consistency of the axiom of choice. (2) It follows from the axiom of replacement that the ordinals necessary to define all sets of integers which can at all be defined in this way will have an upper limit. I doubt that it will be possible to prove that this upper limit is $\omega_1$, as in the case of the constructible sets.~\cite{Godel1946:Remarks-before-the-Princeton-bicentennial-conference-on-problems-in-mathematics}
\end{quote}
Gödel's first question can be seen as fulfilled by the observation that $\HOD\satisfies\ZFC$, which is easy to establish once one knows that $\HOD$ is a definable class. Strictly speaking, however, Gödel referred to the ordinal definable sets, that is, to \OD\ rather than \HOD, and it isn't generally true that \OD\ is a model of \ZF\ set theory. Indeed, \OD\ is not necessarily transitive and need not even satisfy the axiom of extensionality. For example, after forcing to add a Cohen real, not every real number is ordinal definable, and so the set $\R$ of all real numbers and the set $\R^{\HOD}$ of all ordinal-definable real numbers are distinct sets in \OD, but they have the same \OD-members, which violates extensionality. A similar argument with the empty set $\emptyset$ and the set of minimal-rank non-\OD\ sets shows that \OD\ satisfies extensionality if and only if every set is ordinal definable, in which case $V=\HOD=\OD$. 


To summarize the history: 
Gödel introduced the concept of ordinal definability in \cite{Godel1946:Remarks-before-the-Princeton-bicentennial-conference-on-problems-in-mathematics}, but the contemporary approach to \OD\ and \HOD, including the fact that these classes are definable in the language of set theory, appears in~\cite{MyhillScott1971:Ordinal-definability}, who cite independent work of~\cite{Takeuti1961:Remarks-on-Cantors-absolute} and a 1953 preliminary report of Post, as well as unpublished joint work of H\'ajek, \Vopenka, and Balcar.

\section{The forcing modalities and the maximality principles}

The motivating origin of this article was to explore the precise connection between the axiom $V=\HOD$ and the forcing axiom known as the \emph{maximality principle} (\MP), which asserts that every forceably necessary statement is true. That is, if a statement can be forced in such a way that it will remain true in all further forcing extensions, then it is already true. This axiom is naturally expressed via the forcing modalities as the scheme of assertions: 
 $$\possible\necessary\varphi\implies\varphi,$$
where $\possible\psi$ means that $\psi$ is forceable, being forced by some set-sized forcing notion, and $\necessary\psi$ means that $\psi$ is forced by every such forcing notion. In light of the definability of the forcing relation, these modal operators are themselves expressible in the first-order language of set theory. Under the maximality principle the collection of necessary truths is as large as it could be.

The maximality principle is equivalently expressed by the scheme asserting that every true statement is \emph{resurrectible}, that is, it is necessarily possible:
 $$\varphi\implies\necessary\possible\varphi.$$
This is simply the contrapositive of the earlier formulation in the instance with~$\neg\varphi$. The boldface maximality principle $\MP(\R)$ allows real parameters $r$ into the scheme: 
 $$\forall r\in\R\bigl(\possible\necessary\varphi(r)\implies\varphi(r)\bigr).$$
Uncountable parameters can never be allowed into the maximality principle, since it is forceably necessary that any given parameter is countable, and so this would already have to be true.

The maximality principle is easily seen to imply that $V\neq L$, since this is forceable and once true it remains true in all further forcing extensions, so by the maximality principle, it must already be true. More generally, the maximality principle implies that the universe is not a \emph{bedrock}, a model with no nontrivial ground models, since after nontrivial forcing there is a nontrivial ground, of course, and this remains true in all subsequent forcing extensions. So \MP\ implies the failure of the ground axiom \GA, which asserts that the universe was not obtained by nontrivial forcing~\cite{Reitz2006:Dissertation, Reitz2007:TheGroundAxiom}. In particular, under \MP\ the universe is not the \emph{mantle}, $V\neq\M$, defined as the intersection of all ground models---the mantle $\M$ is, provably, the largest forcing-invariant class~\cite{FuchsHamkinsReitz2015:Set-theoreticGeology}. Similarly, the maximality principle implies that $V$ is not equal to many of the other commonly considered canonical inner models, such as the minimal model of a measurable cardinal $V\neq L[\mu]$ or the core model $V\neq K$ and so forth. These models are absolute to any forcing extension, so it is forceably necessary that the universe is not equal to any of them.

What about $V\neq\HOD$? Of course the class $\HOD$ of hereditarily ordinal-definable sets is a centrally important canonical inner model, so in analogy with the results of the previous paragraph one might reasonably expect that \MP\ should imply $V\neq\HOD$. However, a critical difference is that $\HOD$ is not absolute by forcing, and so the simple arguments above showing that \MP\ implies $V\neq L$ and $V\neq \M$ do not seem to show directly that it implies $V\neq\HOD$.

\begin{question}
 What is the precise relationship between the maximality principle \MP\ and the axiom $V=\HOD$?
\end{question}

We shall give a full account, based on the class-forcing contingency of the set-forcing contingency of $V=\HOD$. To be more precise, a statement $\varphi$ is \emph{contingent} by forcing if both $\varphi$ and $\neg\varphi$ can be forced, and the statement is \emph{necessarily} contingent, if this is true in every forcing extension. (Note: In this paper we shall always intend \textbf{set-sized forcing} when speaking of forcing, unless we specifically mention class forcing.) For example, the continuum hypothesis is necessarily contingent by forcing, since both $\CH$ and $\neg\CH$ are forceable over every model of \ZFC. The necessarily contingent statements $\varphi$ are also known as \emph{switches} in the modal logic of forcing since one can turn them successively on and off by forcing, a property that can be expressed via the forcing modalities by:
  $$\necessary\bigl(\possible\varphi\wedge
      \possible\neg\varphi\bigr).$$

We shall prove in \ZFC\ that $V=\HOD$ is a switch whenever it is true, that is, it is necessarily contingent when true. Indeed, it is necessarily contingent whenever it is forceable, and furthermore, it is necessarily contingent whenever it holds somewhere in the generic multiverse, as we shall prove in theorems \ref{Theorem.V=HOD-if-true-is-necessarily-forceable} and \ref{Theorem.TFAE-HOD-is-forceable}. We shall also show in contrast that class forcing can achieve a model of \ZFC\ in which it is not contingent, so $V\neq\HOD$ throughout the generic multiverse; see theorem~\ref{Theorem.VneqHOD-can-be-necessary}. 

\section{A review of the contingency of $V=\HOD$}

Let us now get into the forcing contingency of $V=\HOD$. We begin by reviewing several well-known facts.

\begin{observation}\label{Observation.VneqHOD-is-forceable}
 $V\neq\HOD$ is easily forceable over any model of \ZFC.
\end{observation}

\begin{proof}
This is a well-known observation. One achieves $V\neq\HOD$ after the forcing merely to add a Cohen real, or indeed, after any nontrivial homogeneous forcing. A forcing notion $\P$ is \emph{homogeneous} if for any two conditions $p,q\in\P$ then there are stronger conditions $p'\leq p$ and $q'\leq q$ such that $\P\restrict p'$ is forcing-equivalent to $\P\restrict q'$. This is true of the forcing to add a Cohen real, since all the lower cones of this forcing notion $\Add(\omega,1)\restrict p$ are isomorphic to the whole forcing $\Add(\omega,1)$ and hence also to each other.

The homogeneity property has a key consequence that you can't have a condition forcing a statement $\varphi(\check a)$ about ground-model parameters while another condition forces the negation $\neg\varphi(\check a)$, since the corresponding stronger conditions $p',q'$ are forcing equivalent and therefore will force the same sentences. The key conclusion is that for any assertion $\varphi(\check a)$ involving only check-name parameters, either every condition forces $\varphi(\check a)$ or every condition forces $\neg\varphi(\check a)$.

The further consequence of this for definability is that if $\P$ is homogeneous and $G\of\P$ is $V$-generic, then $\HOD^{V[G]}\of\HOD(\P)^V$. To see this, it suffices by coding sets with sets of ordinals to show that the ordinal definable sets of ordinals in $V[G]$ are in $\HOD(\P)^V$. If $A$ is a set of ordinals defined in $V[G]$ by $\xi\in A\iff\varphi(\xi,\alpha)$, with parameter $\alpha$, then by the key consequence we can define $A$ in $V$ as the set of $\xi$ for which $1\forces_\P\varphi(\check\xi,\check\alpha)$. This defines $A$ in $V$ using parameters $\alpha$ and $\P$.

In particular, the further consequence implies $\HOD^{V[G]}\of V$, so if the forcing is nontrivial and homogeneous, then $V\neq\HOD$ in the forcing extension $V[G]$.
\end{proof}

Next, it is also well known that $V=\HOD$ is forceable by definable class forcing over any model of \ZFC. 

\begin{observation}
 There is a definable proper-class progressively closed forcing notion that forces $V=\HOD$.
\end{observation}

\begin{proof}
The result was first established by McAloon~\cite{McAloon1971:ConsistencyResultsAboutOrdinalDefinability}. Set theorists often used to achieve this by constructing carefully designed forcing iterations that systematically code every set of ordinals, including the new sets added during the course of the iteration, into a feature of the set-theoretic universe, such as into the \GCH\ pattern on the regular cardinals. One would craft a certain bookkeeping function that would tell the forcing iteration exactly which information to code at each stage. The main point was that a properly designed bookkeeping function would ensure that every set added by the iteration eventually has its turn to be coded, which would make it definable in the forcing extension, and thus we get $V=\HOD$ in the extension.

However, one may dispense with the bookkeeping---there is a slick method to ensure that every set is coded, while requiring no intricate advance planning. We simply perform the full-support $\Ord$-iteration $\P$ which at each regular cardinal stage $\delta$ forces with the \emph{lottery sum} of the two coding possibilities---we place the two forcing notions side-by-side and join them with a new $\one$ above. The generic filter in effect chooses one of them generically; see~\cite{Hamkins2000:LotteryPreparation} for a general account of iterated lottery forcing. For example, to use \GCH\ coding this way, we perform an iteration that at regular cardinal stage $\delta$ uses the lottery sum of two forcing notions, one ensuring $2^\delta=\delta^+$ and the other ensuring $2^\delta=\delta^{++}$. The method is flexible and works with the usual alternative coding methods, such as $\Diamond^*_\kappa$ coding. One can also spread out the forcing to occur along an absolutely definable sequence $\delta_\alpha$ of coding points, such as the successors of strong limit cardinals.

The reason this works without any bookkeeping is that the generic filter $G\of\P$ will decide generically whether the \GCH\ holds at $\delta$ or not, by opting for one branch or the other in each lottery. As a consequence, for every set that is added, there is a dense subclass of conditions that make the choices on a later available block of coding points so as to ensure that this set is coded on that block. Basically, we don't need to organize the detailed progression of coding in advance with a bookkeeping function, because the generic filter will ensure that every set is coded somewhere by genericity---it is dense that any given set is coded somewhere, and indeed it will be coded unboundedly often. So this achieves $V=\HOD$ in the class forcing extension.
\end{proof}

The previous argument shows how to force $V=\HOD$ by forcing what is called the \emph{continuum coding axiom} \CCA, which asserts that every set of ordinals appears coded into some interval of the \GCH\ pattern on regular cardinals. Reitz observed that this implies the ground axiom, since every ground $W\of V$ would have to agree with $V$ on the \GCH\ pattern above the size of the forcing realizing $V$ over $W$, and so $W$ would already have every set of ordinals in $V$, which implies $V=W$~\cite{Reitz2006:Dissertation, Reitz2007:TheGroundAxiom}. Another way to express the fact that $\CCA\implies\GA$ is to say that nontrivial forcing always destroys \CCA. In this sense, \CCA\ is extremely fragile whenever it holds. If \CCA\ is true and you perform nontrivial set forcing, you will never get it back by further set forcing. The \CCA\ can hold at most once in the generic multiverse---it can hold only in the bedrock model. Some models of set theory, the ``bottomless'' models, admit no bedrock, and in this case, the \CCA\ fails throughout the generic multiverse.

Let us also discuss a convenient new method of forcing a form of ordinal definability, while preserving all values of the continuum function $\kappa\mapsto 2^\kappa$. Namely, we denote by $\HOD_b$ the class of sets that are hereditarily definable from ordinal parameters and the parameter $b$. We emphasize that with $\HOD_b$ we have just the one parameter $b$, not a set of parameters, so we can still define a version of the global well-order and this will therefore be a \ZFC\ model. It would be equivalent to allow finitely many parameters. In \ZFC\ it is equivalent to say the universe is ordinal-definable from one parameter or from a set of parameters, since we can well order the set of parameters and thereby pick out any one of them with an ordinal. In \ZF, however, the distinction between ordinal-definability with respect to one parameter versus a set of parameters is meaningful. 

\begin{theorem}\label{Theorem.Forcing-V=HOD-while-preserving-continuum-function}
Assume the \SCH\ holds. Then there is a definable, progressively closed class forcing notion with $V[H]=\HOD_b^{V[H]}$ in the forcing extension, using parameter $b=\R^V$, while preserving all cardinals, cofinalities, and also all values of the continuum function $\kappa\mapsto 2^\kappa$. In particular, the forcing preserves the \GCH\ if this holds in the ground model.
\end{theorem}

\begin{proof}
Let $\P$ be the Easton-support forcing $\Ord$-iteration that begins by adding a Cohen real, and thereafter forces at each singular strong limit stage $\delta$ with the lottery sum of $\Add(\delta^+,1)$ and trivial forcing. Let $H\of\P$ be $V$-generic and consider the extension $V[H]$. We may factor this forcing at any stage $\delta$ as the set-sized forcing $\P_\delta$ up to stage $\delta$ followed by the forcing $\P^\delta$ at stage $\delta$ and beyond. Since this tail forcing is $\leq\delta$-closed, this is altogether a progressively closed iteration and therefore preserves \ZFC.

If $\delta$ is a nontrivial stage of forcing, that is, $\delta$ is a singular strong limit, then by the \SCH\ we know $2^\delta=\delta^+$ and consequently the forcing $\Q_\delta$ at stage $\delta$ has size $\delta^+$.  And $\P_\delta$ will have size at most $2^\delta$, which is $\delta^+$. It follows, we claim, that for every regular cardinal $\gamma$, the forcing $\P_\gamma$ has size at most $\gamma$, since either the nontrivial stages of forcing are bounded below $\gamma$, or $\gamma$ is a regular strong limit, hence inaccessible, and in the Easton support we have taken a direct limit of small forcing notions, making $\P_\gamma$ size $\gamma$. 

Let us argue that all regular cardinals are preserved. Suppose that $\kappa$ is a regular cardinal in $V$ and $\gamma<\kappa$ is also regular. Factor the forcing at $\gamma$ as $\P_\gamma*\P^\gamma$. The tail forcing is $\leq\gamma$-closed, so it adds no new $\gamma$ sequences of ordinals and thus in particular, it cannot collapse the cofinality of $\kappa$ to $\gamma$. And the forcing $\P_\gamma$ has size at most $\gamma$, and so it also preserves the regularity of $\kappa$. So $\kappa$ remains regular in $V[H]$. 

A similar argument shows that the values $2^\kappa$ are preserved. Consider any infinite cardinal $\kappa$. The tail forcing at $\kappa$ and beyond is $\leq\kappa$-closed and so does not add new subsets to $\kappa$. And the forcing up to stage $\kappa$ has size altogether at most $2^\kappa$, so it cannot increase the value of $2^\kappa$. So we have preserved the value of $2^\kappa$.

Finally, we argue that $V[H]=\HOD_b^{V[H]}$ in $V[H]$, using parameter $b=\R^V$. Because we added a Cohen real and then performed countably closed forcing, it follows that this forcing has the $\omega_1$-approximation and cover properties, and so the ground model $V$ is a definable class in $V[H]$ using parameter $P(\omega)^V$, which is equidefinable with $b=\R^V$. Next, let us observe that if $V[H]=V[H']$ for $V$-generic filters $H,H'\of\P$, then $V[H_\gamma]=V[H'_\gamma]$ for every stage of forcing $\gamma$ along the way, we claim, where $H_\gamma\of\P_\gamma$ refers to the forcing iteration prior to $\gamma$. To see this, we may assume inductively that this is true below $\gamma$. If there is a largest stage of nontrivial forcing below $\gamma$, then we get the conclusion immediately. So we may assume that $\gamma$ is a limit of nontrivial stages of forcing, so in particular, $\gamma$ itself will be a strong limit cardinal. Since $H_\gamma\in V[H]=V[H']$ and it is generated from the conditions in the prior sets $H_\beta$ for $\beta<\gamma$, which form a set of size $2^{<\gamma}=\gamma$, it follows that $H_\gamma$ must be an element of $V[H'_\gamma]$, since it can't be added by the later forcing. So $V[H_\gamma]\of V[H'_\gamma]$. The converse holds similarly, so we have $V[H_\gamma]=V[H'_\gamma]$, as desired. A similar argument shows that if $V[H_\theta]=V[H'_\theta]$ for filters $H_\theta,H'_\theta\of\P_\theta$ performing the forcing up to stage $\theta$, then $V[H_\gamma]=V[H'_\gamma]$ for each $\gamma\leq\theta$. 
It follows that $H$ opts for the nontrivial forcing at stage $\delta$, that is, adding a subset to $\delta^+$, if and only if $H'$ also opted this way, since otherwise we couldn't have $V[H_{\delta+1}]=V[H'_{\delta+1}]$. 

Therefore, we claim, the class $A\of\Ord$ consisting of the stages $\delta$ at which the generic filter opted for the nontrivial forcing in the stage $\delta$ lottery is definable from $b$ in $V[H]$. Namely, from $b$ we can define $V$, and then observe that $\delta\in A$ if and only if there is some $H'_{\delta+1}\of\Q_\delta$ that is $V$-generic and opted for nontrivial forcing in the stage $\delta$ lottery and $V[H'_{\delta+1}]$ contains the full $P(\delta)^{V[H]}$. (Note that we are not using $H$ as a parameter here, but rather evaluating truth in $V[H]$, which is the ambient universe in which the definition takes place.) The point is that the argument of the previous paragraph shows that all generic filters have to agree on whether they opted for nontrivial forcing in the lottery at stage $\delta$. They might use different generic filters at that stage, but either they all opt for triviality or they all opt for nontriviality in that lottery.

Finally, we observe that the definability of $A$ from $b$ in $V[H]$ causes $V[H]=\HOD_b^{V[H]}$, since any set of ordinals that is added at any stage will have a chance to be coded into the pattern of membership in $A$ of the singular strong limit cardinals. Basically, $A$ is a $V$-generic Cohen class of ordinals, which therefore encodes every set of ordinals in its diverse intervals. So $V[H]=\HOD_b^{V[H]}$, as desired.
\end{proof}

It follows from the main theorem that with further set-sized forcing, we can achieve the axiom $V=\HOD$, although this may not preserve all values of $2^\kappa$. 

\section{The axiom $V=\HOD$ can be impossible by set forcing}\label{Section.V=HOD-can-be-impossible}

We have just described a few ways to force $V=\HOD$ with class forcing. But what about set forcing? Since we aimed in the class iteration to code every set in the universe, the proper-class nature of the forcing may seem inherent to achieving this. Can we expect ever to force $V=\HOD$ by set-sized forcing, if the axiom does not already hold in the ground model? For example, if $V=\HOD$ holds and we add a Cohen real, then we have $V\neq\HOD$, but can we recover $V=\HOD$ with set forcing? We have already observed that we cannot recover \CCA\ by further set forcing and so we cannot expect to force $V=\HOD$ that way. Is there another way?\label{Page.Forcing-V=HOD-by-set-forcing}

We will answer that question in the next section, but in this section, let us show that the class-forcing aspect of forcing $V=\HOD$ is sometimes required. Every model of \ZFC, we claim, has a class forcing extension in which $V\neq\HOD$ holds in all subsequent set forcing extensions.

\begin{theorem}\label{Theorem.VneqHOD-can-be-necessary}
 The set-theoretic universe $V$ has a class-forcing extension $V[G]$ such that $V=\HOD$ is not forceable over $V[G]$ by set-sized forcing. In short:
  $$V[G]\satisfies\necessary (V\neq \HOD).$$
\end{theorem}

\begin{proof}
We start in $V\satisfies\ZFC$. Let $\P$ be the Easton-support product forcing that adds a Cohen subset to every infinite regular cardinal $\P=\prod_\delta \Q_\delta$, where $\Q_\delta=\Add(\delta,1)^V$. Let $G\of\P$ be $V$-generic and argue in $V[G]$. Since this is homogeneous forcing, we have $V\neq\HOD$ in $V[G]$. 

We claim that every forcing extension of $V[G]$ by set-sized forcing will also have $V\neq\HOD$. Suppose that $V[G][h]$ is a forcing extension via $h\of\Q\in V[G]$. Choose a regular cardinal $\delta$ much larger than $\Q$ and let $g_\delta\of\Add(\delta,1)$ be the generic added by $G$ on coordinate $\delta$. We can rearrange the full forcing extension as 
 $$V[G][h]=V[G^*_\delta][h][g_\delta],$$
where $G^*_\delta$ is the same as $G$, except that the forcing at coordinate $\delta$ is omitted. We have then reinserted this stage of forcing separately at the end with $g_\delta$. The forcing $\Q_\delta=\Add(\delta,1)^V$ was homogeneous in $V$. Although this forcing is no longer Cohen forcing after the $h$ forcing, nevertheless, it remains homogeneous in $V[G^*_\delta][h]$, and so we shall have $V\neq\HOD$ in the full extension $V[G^*_\delta][h][g_\delta]=V[G][h]$. So $V\neq\HOD$ in every set-forcing extension of $V[G]$. 
\end{proof}

\section{When true, $V=\HOD$ is a switch}

We claim next, in contrast, that it is also possible to have a model of set theory for which $V=\HOD$ is forceable by set-sized forcing over every forcing extension. In fact, we claim, whenever $V=\HOD$ holds, then in every forcing extension we can recover $V=\HOD$ with further set-sized forcing. In particular, after adding a Cohen real over a model of $V=\HOD$, we can indeed recover $V=\HOD$ with set forcing, which therefore answers the question we had posed about this at the beginning of section~\ref{Section.V=HOD-can-be-impossible}.\goodbreak

\begin{theorem}\label{Theorem.V=HOD-if-true-is-necessarily-forceable}
 $V=\HOD$, if true, is forceable over every forcing extension. 
 $$V=\HOD\quad\Implies\quad\necessary\possible (V=\HOD).$$ 
Thus the axiom $V=\HOD$, when true, is a switch.
\end{theorem}

\begin{proof}
This observation was made by the first author in \cite{Hamkins2015:HOD-of-a-set}. Suppose that $V=\HOD$ holds in $V$. We want to show that it remains forceable over any forcing extension $V[g]$, where $g\of\P\in V$ is $V$-generic. For this, we shall perform \emph{self-encoding forcing} over $V[g]$ above the size of $\P$; see~\cite{FuchsHamkinsReitz2015:Set-theoreticGeology} for an introduction to self-encoding forcing. This will provide a subsequent forcing extension $V[g][h]$ in which $g$ and $h$ become ordinal definable in the extension. And we can also arrange a bit more definability, which will imply $V=\HOD$ in the extension $V[g][h]$.

To explain the details, we pick a $\beth$-fixed point $\lambda_0$ with $\P\in V_{\lambda_0}$ and then design a forcing notion $\Q_0$ that will ensure that $V_{\lambda_0}[g]$ is coded into the \GCH\ pattern at the successive regular cardinals above $\lambda_0$ in the resulting forcing extension $V[g][h_0]$. Then, let $\lambda_1$ be first $\beth$-fixed point above $\Q_0$ and force to add $h_1$ that codes $h_0$ into the next block of cardinals. Then we add $h_2$ to code $h_1$ in a higher block of cardinals, and so forth for the full-support $\omega$-iteration $\Q$, adding $h_{n+1}$ so as to code $h_n$ into the \GCH\ pattern above $\lambda_{n+1}$. The iteration $\Q$ forces to add a filter $h\of\Q$ over $V[g]$ such that in $V[g][h]$, the sets $V_{\lambda_0}$, $g$, and also $h$ itself are coded into the \GCH\ pattern. Also $\P$ becomes ordinal definable, being in $V_{\lambda_0}$, and $\Q$ is ordinal definable in the extension as well.

Since the forcing $\Q$ is $\leq\lambda_0$-closed, it follows that $\P*\Qdot$ admits a closure point at $\lambda_0$, in the sense of~\cite{Hamkins2003:ExtensionsWithApproximationAndCoverProperties}. This extension therefore admits the $\lambda_0$-approximation and cover properties and so by the ground-model definability theorem (see~\cite{Hamkins2003:ExtensionsWithApproximationAndCoverProperties, Laver2007:CertainVeryLargeCardinalsNotCreated, FuchsHamkinsReitz2015:Set-theoreticGeology}) we may conclude that $V$ is definable in $V[g][h]$ from parameter $V_{\lambda_0}$. But since this set $V_{\lambda_0}$ is ordinal definable in $V[g][h]$, it follows that the ground model $V$ is an ordinal-definable class in $V[g][h]$. And since $V$ satisfies $V=\HOD$, it therefore follows that $V\of\HOD^{V[g][h]}$. And since $g\of\P$ and $h\of\Q$ also are ordinal-definable in $V[g][h]$, it follows that 
$V[g][h]=\HOD^{V[g][h]}$, as desired. 

For the final remark, $V=\HOD$ will be a switch whenever it is true, since the previous argument shows how to turn it back on over any forcing extension, and by observation \ref{Observation.VneqHOD-is-forceable} it is easily turned off by any nontrivial homogeneous forcing.
\end{proof}

\section{The contingent \HOD\ dichotomy}\label{Section.The-contingent-HOD-dichotomy}

We should like now to sharpen the analysis by setting it in the context of the \emph{generic multiverse} of $V$, the collection of set-theoretic universes reachable in finitely many steps of passing successively to a forcing extension or to a ground. By the downward directed grounds hypothesis, established by Usuba~\cite{Usuba2017:The-downward-directed-grounds-hypothesis-and-very-large-cardinals}, we know that the grounds of any model of set theory are downward directed, and by iterative application of this it follows that every model in the generic multiverse of $V$ is a forcing extension $W[g]$ of a ground $W$ of $V$. Thus, to exhaust the generic multiverse, two steps suffice: go down to a suitable ground $W$ and then up to a suitable forcing extension $W[g]$.

Part of our analysis will involve the \emph{blurry} \HOD\ hierarchy, introduced by Gunter Fuchs \cite{Fuchs2022:Blurry-definability, Fuchs2024:On-the-consistency-strength-of-critical-leaps, Fuchs2025:More-on-blurry-HOD}, which consists of the tower of inner models $\kappa\blurryHOD$. The core idea generalizes the concept of ordinal-algebraicity, introduced by the first author and Cole Leahy \cite{HamkinsLeahy2016:AlgebraicityAndImplicitDefinabilityInSetTheory}, by which a set is \emph{ordinal algebraic} when it is a member of a finite ordinal-definable set, and $\HOA$ is the class of hereditarily ordinal-algebraic sets---the central theorem establishes $\HOA=\HOD$. More generally, we say that a set is \emph{$\kappa$-algebraic} when it is a member of an ordinal definable set of size less than $\kappa$, and $\kappa\blurryHOD$ consists of the sets such that they and all their hereditary members are $\kappa$-algebraic. Fuchs shows that $\kappa\blurryHOD$ is an inner model of \ZF, though not necessarily \ZFC. It is easy to see in \ZFC\ that the blurry hierarchy accumulates to the whole universe $V=\Union_\kappa \kappa\blurryHOD$, and this hierarchy is said to \emph{collapse} if $V=\kappa\blurryHOD$ for some cardinal $\kappa$. 

\begin{theorem}\label{Theorem.TFAE-HOD-is-forceable}
The following are equivalent in \ZFC.
 \begin{enumerate}
  \item The inner model $\HOD$ is a ground of $V$.
    $$V=\HOD[G]$$
  \item The axiom $V=\HOD$ is forceable.
     $$\possible(V=\HOD)$$
  \item The universe is ordinal-definable from a set. $$\exists b\ V=\HOD_b$$ 
  \item The universe admits a definable global well order (allowing parameters).
  \item The blurry \HOD\ hierarchy collapses. 
   $$\exists \kappa\ V=\kappa\blurryHOD$$
  \item The axiom $V=\HOD$ holds somewhere in the generic multiverse. 
  \item For every model $W$ in the generic multiverse, $\HOD^W$ is a ground.
   $$W=\HOD^W[G]$$
  \item The axiom $V=\HOD$ is forceable over every $W$ in the generic multiverse.
  \item Every model $W$ in the generic multiverse is ordinal-definable from a set.
   $$\exists b\ W=\HOD_b^W$$
  \item Every model in the generic multiverse admits a definable global well order definable (allowing parameters).
  \item In every model $W$ of the generic multiverse, the blurry \HOD\ hierarchy collapses. 
    $$\exists\kappa\ W=\kappa\blurryHOD^W$$
 \end{enumerate}
\end{theorem}

\begin{proof}
The equivalences of 2, 3, 6, 8, and 9 are proved in \cite{Hamkins2015:HOD-of-a-set}; the connection with 1 (and consequently also 7) is provided by \cite{Woodin_Davis_Rodriguez_2012}, credited to \Vopenka; and the connection with 5 and hence also 11 is provided by \cite{Fuchs2022:Blurry-definability}. The equivalences $(3\iff 4)$ and $(9\iff 10)$ are well known. We shall make a completely self-contained proof of everything here.

$(1\implies 3)$ If $\HOD$ is a ground, then $V=\HOD[G]$ for some $\HOD$-generic $G\of\P\in\HOD$, so clearly $V=\HOD_G$. 

$(3\iff 4)$ This is well known. If $V=\HOD_b$, then we can define the parameterized \HOD\ order, by which $x\leq y$ just in case $x$ is ordinal definable from $b$ in some smaller $V_\theta$ than $y$ is, or they use the same $V_\theta$ but $x$ uses a smaller defining formula, or they use the same $V_\theta$ and the same formula, but $x$ uses  smaller ordinal parameters. Conversely, if there is a definable well ordering using parameter $b$, then by also stratifying by rank there is one in order type $\Ord$, so every object is the $\alpha$th object in this order for some ordinal $\alpha$, and so $V=\HOD_b$. Note that the order we defined assuming 3 had complexity $\Sigma_2(b)$, and so 4 is equivalent to asserting that there is a $\Sigma_2(b)$-definable well ordering of the universe, which is therefore expressible in the language of set theory.

$(3\implies 2)$ We adapt the proof of theorem \ref{Theorem.V=HOD-if-true-is-necessarily-forceable}. Suppose that $V=\HOD_b$ for some set $b$. We may assume without loss that $b$ is a set of ordinals, since such sets can code any given set. We use the method of self-encoding forcing, as described in the proof of theorem \ref{Theorem.V=HOD-if-true-is-necessarily-forceable}. Specifically, we consider the forcing $h_0\of\Q_0$ that first adds a Cohen real and then performs forcing that codes $P(\omega)^V$ and the parameter $b$ into the \GCH\ pattern on regular cardinals starting at a large beth fixed point $\lambda_0$, adding generic filter $h_0$; and then performs forcing $h_1\of\Q_1$ to code $h_0$ into the \GCH\ pattern starting at the next beth fixed point $\lambda_1$; and so forth adding $h_{n+1}$ to code $h_n$. The full-support $\omega$-iteration of these forcing notions has a closure point at $\omega$, since it added a Cohen real followed by countably closed forcing, and so $V$ is definable in $V[h]$ from the parameter $P(\omega)^V$, which itself is definable in $V[h]$. Since $b$ is also definable, we see that every element of $V$ becomes ordinal definable in $V[h]$. And since $h$ is definable as well, we conclude that $V[h]$ satisfies $V=\HOD$.

$(2\implies 3)$ Suppose that $V=\HOD$ holds in some forcing extension $V[g]$ with $V$-generic $g\of\Q\in V$. Let $b\in V$ be a set of ordinals coding $\Q$ and every condition $p\in\Q$. Consider any set of ordinals $A\of\Ord$ in $V$. It becomes ordinal-definable in $V[g]$, so $V[g]\satisfies\alpha\in A\iff\varphi(\alpha,\theta)$ for some assertion $\varphi$ with ordinal parameter $\theta$. This situation is forced by some condition $p$, and so we may observe that $\alpha\in A\iff p\forces\varphi(\check\alpha,\check\theta)$, which is definable in $V$ from parameters $b$ and $\theta$ and the ordinals required to define $p$ from $b$. So every set of ordinals in $V$ is ordinal definable from $b$, and since every set is coded by a set of ordinals, this means $V=\HOD_b$. 

$(3\implies 5)$ Assume that $V=\HOD_b$ for some set $b$. We now simply ``blur'' the parameter set, namely, let $B$ be the set of minimal rank sets $b$ for which $V=\HOD_b$, which is a definable set, and let $\kappa>|B|$. Given any set $a$, since it is in $V=\HOD_b$, we know $a$ is defined by the property $\varphi(a,b,\vec\alpha)$ for some ordinal parameters. Let $A$ be the set of all $a'$ that are defined by $\varphi(a',b',\vec\alpha)$ for some $b'\in B$. This is an ordinal-definable set of size less than $\kappa$, so $a\in\kappa\blurryHOD$, as desired.

$(5\implies 3)$ Although not needed for the cycle of implications in the theorem here, let us give a direct proof of this converse implication, just to illustrate the ideas further. Assume \ZFC\ and $V=\kappa\blurryHOD$ for some cardinal $\kappa$. Let $b$ be a well-ordered enumeration of $2^{<\kappa}$. We claim that $V=\HOD_b$. To see this, it suffices by \AC\ to show that every set of ordinals $a$ is in $\HOD_b$. Suppose $a\of\lambda$ for some ordinal $\lambda$. By our assumption $a$ is $\kappa$-algebraic, so $a\in A$ for some ordinal-definable set $A$ of size less than $\kappa$. We may assume that $A$ consists of subsets of $\lambda$. For each distinct pair $x,y\in A$, let $\alpha_{x,y}=\min x\triangle y$, the smallest ordinal distinguishing $x$ and $y$, and place these all into the set $I=\set{\alpha_{x,y}\mid x,y\in A, x\neq y}$; a similar method is used in \cite[theorem~3]{Fuchs2022:Blurry-definability}. Note that each $x\in A$ is determined uniquely by its pattern on $I$, that is, by $x\intersect I$, since we distinguish $x$ from every $y\neq x$ at $\alpha_{x,y}$. Furthermore, the set $I$ is definable from $A$, which is itself ordinal definable, and so the collection $A_0=\set{ a \intersect I\mid a\in A}$ is ordinal definable. This set $A_0$ is definably isomorphic to a subset of $P(\gamma)$, where $\gamma$ is the order type of $I$, which is an ordinal less than $\kappa$. The parameter $b$ provides a well-order of that set, and so we may definably pick out $a\intersect I$ with the appropriate ordinal index. Since $a\intersect I$ definably determines $a$, it follows that $a$ is ordinal definable with parameter $b$. So $V=\HOD_b$, as desired.

$(5\implies 1)$ This is essentially \cite[theorem~5]{Fuchs2022:Blurry-definability}, but let us give the argument. We shall show in \ZFC\ that if $V=\kappa\blurryHOD$, then $V$ is a $\kappa$-c.c.~forcing extension of \HOD. We shall do this by verifying \Bukovsky's criterion \cite{Bukovsky1973:Characterizatoin-of-generic-extension-of-models-of-set-theory}, which asserts that $V$ is a $\kappa$-c.c.~forcing extension of inner model $M$ of \ZFC\ precisely when every function $f:\gamma\to\gamma$ in $V$ on an ordinal $\gamma$ is uniformly $\kappa$-covered by a function $F:\gamma\to M$ in $M$, meaning that $F\in M$, that $f(\alpha)\in F(\alpha)$, and that $|F(\alpha)|<\kappa$ for all $\alpha<\gamma$. To verify this with $M=\HOD$, suppose $f:\gamma\to\gamma$. Since $V=\kappa\blurryHOD$, there is an ordinal-definable family $\mathcal F$ of size less than $\kappa$, with $f\in\mathcal F$. Let $F(\alpha)=\set{h(\alpha)\mid h\in\mathcal F}$, which has size less than $\kappa$, and contains $f(\alpha)$. So we have fulfilled \Bukovsky's criterion, and therefore $V$ is a $\kappa$-c.c.~forcing extension of \HOD, as desired. 

$(3\implies 1)$ Although not needed for the cycle of implications to prove the theorem here, let us mention nevertheless that the implication $3\implies 1$ is provided by \cite[theorem~6]{Woodin_Davis_Rodriguez_2012}, credited to \Vopenka.

$(1\iff 2\iff 3\iff 4\iff 5\iff 6\iff 7\iff 8\iff 9\iff 10\iff 11)$ We have already proved the equivalence of 1, 2, 3, 4, and 5. So it will suffice to show that these statements are invariant by forcing, preserved both to forcing extensions and grounds, for in this case, if 1, 2, 3, 4, or 5 hold anywhere in the generic multiverse, then they will all propagate to hold throughout the generic multiverse, which is what the other statements 6--11 express. Statement 2 is obviously preserved to grounds. And statement 3 is preserved to forcing extensions. To see this, suppose $V=\HOD_b$ and consider forcing extension $V[g]$ via generic filter $g\of\Q\in V$. By the ground-model definability theorem, $V$ is definable in $V[g]$ from parameter $P(\Q)^V$, so with parameters $b$ and $P(\Q)^V$ we can define $\Q$ and the class of names $V^\Q$, and so with ordinals as well we can define any particular $\Q$-name in $V$, and so with parameter $g$ in addition we can define every particular element of $V[g]$. So 3 holds in $V[g]$. Alternatively, we could have argued that 2 is preserved by forcing, since if $V=\HOD$ is forced by some forcing notion $\Q$, then for any forcing extension $V[g]$, we may perform further forcing with $h\of\Q$ to form $V[g][h]$, which is an extension of $V[h]$, where $V=\HOD$ holds. By theorem \ref{Theorem.V=HOD-if-true-is-necessarily-forceable}, the axiom $V=\HOD$ is forceable over every forcing extension of $V[h]$, and so it is forceable over $V[h][g]$ and hence also over $V[g]$. So we have proved that statements 1, 2, 3, 4, and 5 are preserved to forcing extensions and grounds, and thus hold throughout the generic multiverse, if they hold at all, and this gives the equivalence with statements 6 through 11.
\end{proof}

This analysis reveals a fundamental dichotomy concerning the nature of $\HOD$ in the generic multiverse, which we should like to highlight. We call it the \emph{contingent \HOD\ dichotomy}:

\begin{corollary}
 \ZFC\ proves the following dichotomy. Either
 \begin{enumerate}
  \item The following statements hold in $V$ and throughout the generic multiverse.
  \begin{enumerate}
    \item The inner model $\HOD$ is a ground.
    \item $\possible(V=\HOD)$.
    \item The universe is ordinal-definable from a set.
    \item There is a definable global well order (from parameters).
    \item The blurry \HOD\ hierarchy collapses, so $V=\kappa\blurryHOD$ for some $\kappa$. 
  \end{enumerate}
\smallskip
\hskip-3.5em Or
  \item \vskip-\baselineskip
  Those statements fail in $V$ and throughout the generic multiverse. 
  \begin{enumerate}
    \item The inner model \HOD\ is not a ground.
    \item $\necessary(V\neq\HOD)$.
    \item The universe is not ordinal-definable from a set. 
    \item There is no definable global well order (from parameters).
    \item The blurry \HOD\ hierarchy doesn't collapse, so $V\neq\kappa\blurryHOD$ for all $\kappa$.
  \end{enumerate} 
 \end{enumerate}
If \ZFC\ is consistent, then both prongs of this dichotomy occur.
\end{corollary}

In the first prong of the dichotomy, for example, the universe is realized as a forcing extension $V=\HOD[G]$; there is a forcing extension $V[H]$ satisfying the axiom $V=\HOD$; the universe is definable from a set $V=\HOD_b$ for some parameter $b$; there is a definable well ordering of the universe from some parameter $b$; and the universe is $\kappa\blurryHOD$ for some cardinal $\kappa$; and what is more, these statements also hold of every world $W$ throughout the generic multiverse, so $W=\HOD^W[G]$; some forcing extension $W[G]$ satisfies $V=\HOD$; the universe is $W=\HOD_b^W$; there is a parametrically definable global well order in $W$; and $W=\kappa\blurryHOD^W$ for some cardinal $\kappa$. Meanwhile, in the second prong of the dichotomy, we have $V\neq\HOD$ throughout the generic multiverse; $\HOD^W$ is never a ground of $W$; the universe $W$ is never ordinal-definable from a set; there is never an ordinal-definable global well order; and $W\neq\kappa\blurryHOD^W$ for all $\kappa$. 

\begin{proof}
The dichotomy is the main content of theorem \ref{Theorem.TFAE-HOD-is-forceable}, which shows that in any model of \ZFC\ either all the statements of theorem \ref{Theorem.TFAE-HOD-is-forceable} hold or they all fail. If \ZFC\ is consistent, then the first prong of the dichotomy occurs in any model of $V=\HOD$, and the second occurs in some models in light of theorem \ref{Theorem.VneqHOD-can-be-necessary}.
\end{proof}

For brevity, we denote the two prongs of the contingent \HOD\ dichotomy by:
 $$\possible(V=\HOD)\qquad\text{ versus }\qquad\necessary(V\neq\HOD)$$
The first prong asserts that the set-theoretic universe $V$ is close to \HOD\ by forcing, in that $V=\HOD$ is forceable, hence contingent, and this is equivalent to \HOD\ being a ground, and these properties are invariant by forcing, true in every world throughout the generic multiverse. The second prong of the dichotomy, in contrast, asserts that $V$ is far from \HOD, that $V=HOD$ is incontingent by forcing, being necessarily false, and this also is invariant by forcing, holding throughout the generic multiverse. 

\section{The contingent \HOD\ dichotomy in canonical inner models}
\label{Section.Contingent-HOD-dichotomy-in-inner-models}

The various canonical inner models commonly studied in set theory generally appear in the first prong of the contingent \HOD\ dichotomy. For example, if the universe is constructible from a set, then it is in the first prong of the dichotomy.

\begin{theorem}\label{Theorem.L[A]-is-in-first-prong}
 If the universe is constructible from a set $V=L[A]$, then it is in the first prong $\possible(V=\HOD)$ of the contingent \HOD\ dichotomy.
\end{theorem}

\begin{proof}
If $V=L[A]$, then $V=\HOD_A$ is true, so we are in the first prong of the dichotomy.
\end{proof}

This already covers many of the canonical inner models of large cardinals, such as $L$, $L_\mu$, and more. But also, of course whenever an inner model satisfies $V=\HOD$, then it will be in the first prong of the dichotomy, which covers many further cases, including many instances of $L[\vec E]$. 

More generally, let us call attention to Gabriel Goldberg's result showing how his ultrapower axiom with a supercompact cardinal implies that we are necessarily in the first prong of the contingent \HOD\ dichotomy. 

\begin{theorem}[{\cite[Theorem~6.2.8]{Goldberg2022UltrapowerAxiom}}]\label{Theorem.Goldberg-ultrapower-axiom}
Assume the Ultrapower Axiom and that there is a supercompact cardinal. Then \HOD\ is a ground of $V$.
\end{theorem}

So these models are also therefore in the first prong of the contingent \HOD\ dichotomy. Goldberg says that ``[the ultrapower axiom] holds in all canonical inner models, a result that is philosophically central to this dissertation,'' and this theorem fulfills this philosophical perspective by showing that under the ultrapower axiom (in the presence of a supercompact cardinal), we are in the close-to-\HOD\ prong of the contingent \HOD\ dichotomy.

\section{The contingent \HOD\ dichotomy and large cardinals}

Next, we prove that all the usual large cardinals are relatively consistent with both sides of the contingent \HOD\ dichotomy. In short, large cardinals do not settle the contingent \HOD\ dichotomy. 

\begin{theorem}\label{Theorem.Large-cardinals-consistent-with-both-prongs}
 All the usual large cardinal hypotheses are consistent with both prongs of the contingent \HOD\ dichotomy.
\end{theorem}

\begin{proof}
It is well known how to force $V=\HOD$ with class forcing while preserving any of the usual large cardinal notions. Of course this is true for large cardinals such as measurable, weakly compact, Ramsey, superstrong, rank-to-rank, and related notions, since these are witnessed in each case inside some large-enough $V_\theta$ (so these notions are $\Sigma_2$-expressible), and we can force $V=\HOD$ in the usual way, but simply starting above $\theta$. But in fact it is also known how to force $V=\HOD$ while preserving a proper class of these notions, as well as strong cardinals, supercompact cardinals, and so forth. This is a staple of the Easton-support master condition lifting arguments of large cardinal forcing theory.

For the other prong, we point out that all the standard large cardinal hypotheses are also preserved by the forcing in theorem \ref{Theorem.VneqHOD-can-be-necessary}, which is an $\Ord$-product adding a Cohen set at regular cardinals (or on a certain class of regular cardinals). So we can also have $\necessary(V\neq\HOD)$ with large cardinals.
\end{proof}

Meanwhile, it may be important to mention that large cardinals need not generally exhibit their large cardinal nature in \HOD. For example, \cite{ChengFriedmanHamkins2015:LargeCardinalsNeedNotBeLargeInHOD} shows that a supercompact cardinal need not be weakly compact in \HOD, and there can be a proper class of supercompact cardinals in $V$, none of them weakly compact in \HOD, with no supercompact cardinals in \HOD. Similar results hold for other large cardinals, such as measurable and strong cardinals.

\section{Contingent \HOD\ dichotomy versus \GCH\ and Easton's theorem}

We know that the contingent \HOD\ dichotomy is invariant by forcing---you can't change which prong you are in by set forcing. Meanwhile, with class forcing, you can change things, since of course you can force $V=\HOD$ and also $\necessary(V\neq\HOD)$ by class forcing. What about the other natural class forcing possibilities? Let us investigate the situation with \GCH.

\begin{theorem}
 The \GCH\ is relatively consistent with both prongs of the contingent \HOD\ dichotomy.
\end{theorem}

\begin{proof}
The forcing of theorem \ref{Theorem.VneqHOD-can-be-necessary} will force the \GCH, since when one adds a Cohen subset to every regular cardinal then you systematically force the \GCH\ at every cardinal along the way. Therefore, $\necessary(V\neq\HOD)$ is relatively consistent with \GCH. On the other hand, we can also have $V=\HOD+\GCH$, since this is true in $L$ and other canonical models, so the first prong is possible with \GCH. 
\end{proof}

More generally, let us observe that Easton's theorem is compatible with either prong of the contingent \HOD\ dichotomy. 

\begin{theorem}
 Assume $\ZFC+\GCH$ and suppose that $E$ is a class function defined on the regular cardinals in accordance with the requirements of Easton's theorem. Then there are cardinal-preserving class forcing extensions $V[G]$ and $V[H]$, both satisfying \ZFC\ and also $2^\kappa=E(\kappa)$ for every regular cardinal, such that $V[G]$ is in one prong of the contingent \HOD\ dichotomy and $V[H]$ is in the other. 
\end{theorem}

\begin{proof}
For the second prong, let us take $V[G]$ to arise from Easton's original forcing, the Easton-support class product $\prod_\kappa\Add(\kappa,E(\kappa))$. This forcing, of course, achieves Easton's result in which $2^\kappa=E(\kappa)$ in the forcing extension $V[G]$ for every regular cardinal $\kappa$. But furthermore, since Easton's forcing is homogeneous, we will have $\HOD^{V[G]}\of V$, and this implies that \HOD\ is not a ground in $V[G]$, since $V$ is not a set-forcing ground of $V[G]$ and therefore can contain no set-forcing ground of $V[G]$. So $V[G]$ must be in the second prong of the contingent \HOD\ dichotomy, with $\necessary(V\neq\HOD)$.

For the first prong, we simply apply theorem \ref{Theorem.Forcing-V=HOD-while-preserving-continuum-function} over the model $V[G]$. Easton's forcing preserves the \SCH, and so by theorem \ref{Theorem.Forcing-V=HOD-while-preserving-continuum-function} there is a further extension $V[G][H]$ forcing $V=\HOD_b$, which places us into the first prong of the contingent \HOD\ dichotomy, while also preserving the values of the continuum function $2^\kappa=E(\kappa)$ that we had so carefully arranged in $V[G]$, as desired.
\end{proof}

It seems to be left open by this particular method whether one can achieve $V=\HOD$ and the Easton result simultaneously over any model of $\ZFC+\GCH$. We expect that this will be possible with more intricate coding methods.

\section{The contingent \HOD\ dichotomy and set-theoretic geology}

Allow us briefly to recall some definitions from set-theoretic geology \cite{FuchsHamkinsReitz2015:Set-theoreticGeology}. A \emph{ground} of $V$ is an inner model $W\satisfies\ZFC$ over which the universe was obtained by forcing $V=W[G]$, for some $W$-generic $G\of\Q\in W$. The \emph{bedrock} is a minimal ground model, which is unique when it exists and furthermore is the least ground, contained in all others \cite{Usuba2017:The-downward-directed-grounds-hypothesis-and-very-large-cardinals}. The \emph{mantle} $\Mantle$ is the inner model resulting from the intersection of all grounds
 $$\Mantle=\Intersect_{V=W[G]} W.$$
This is definable in \ZFC\ in the light of the ground-model enumeration theorem \cite[theorem~12]{FuchsHamkinsReitz2015:Set-theoreticGeology}, and it is a model of \ZFC, a consequence of \cite[theorem~22]{FuchsHamkinsReitz2015:Set-theoreticGeology} and \cite{Usuba2017:The-downward-directed-grounds-hypothesis-and-very-large-cardinals}, but it is not necessarily itself a ground, as some models of \ZFC\ have no minimal ground, causing the mantle to dig deeper and deeper. The \emph{generic \HOD}, introduced by Gunter Fuchs \cite[p.~298]{Fuchs2008:ClosedMaximalityPrinciples}, is the intersection of all the \HOD s arising in any forcing extension:
 $$\gHOD=\Intersect_{V\of V[G]} \HOD^{V[G]}.$$
Since every forcing notion is absorbed by sufficiently large collapse forcing, which is homogeneous, it follows for this further extension that $\HOD^{V[G]}$ is contained in $V$. Somewhat paradoxically, therefore, as one proceeds to further collapse forcing extensions, the \HOD s of those models become smaller. The sets in the generic \HOD\ are precisely those sets that are necessarily ordinal-definable in all forcing extensions; meanwhile, the sets outside \gHOD\ have necessarily contingent ordinal definability---they can become definable in a forcing extension, of course, as any set can, but they will inevitably lose this feature with further forcing. The generic \HOD\ is a definable class in $V$ and furthermore a model of \ZFC\ by \cite{Fuchs2008:ClosedMaximalityPrinciples}; see also \cite[theorem~63]{FuchsHamkinsReitz2015:Set-theoreticGeology}. The results of \cite[theorem~66,67,68,70,71]{FuchsHamkinsReitz2015:Set-theoreticGeology} exhibit various possible behaviors in models of set theory, with $\gHOD=\HOD=\Mantle\ofneq V$, or $\gHOD=\Mantle\ofneq \HOD=V$, or $\gHOD=\HOD\ofneq\Mantle=V$, or 
$\gHOD=\HOD=\Mantle=V$ all occurring in various models of set theory. 

We observe that when there is a bedrock model, a minimal ground model, then the question of which prong of the contingent \HOD\ dichotomy holds is reflected exactly in the question of whether $V=\HOD$ holds in the bedrock.

\begin{theorem}\label{Theorem.Contingent-HOD-dichotomy-with-bedrock}
 If the set-theoretic universe $V$ admits a bedrock, then we are in the $\possible(V=\HOD)$ prong of the contingent \HOD\ dichotomy if and only if the bedrock is a model of $V=\HOD$. If the bedrock satisfies $V\neq\HOD$, then $V\neq\HOD$ holds throughout the generic multiverse.
\end{theorem}

\begin{proof}
If the bedrock satisfies $V=\HOD$, then we are in the $\possible(V=\HOD)$ case of the dichotomy, of course, since the bedrock is one of the models in the generic multiverse. Alternatively, if $V\neq\HOD$ in the bedrock $W$, then in particular, $\HOD^W$ is not a ground of $W$ since the only ground of a bedrock is the bedrock itself, and so we must be in the other prong of the dichotomy.
\end{proof}

Meanwhile, the dichotomy does not settle the question of whether the universe has a bedrock. 

\begin{theorem}
 If \ZFC\ is consistent, then both prongs of the contingent \HOD\ dichotomy are consistent both with the existence and with the nonexistence of a bedrock. All four combinations occur.
\end{theorem}

\begin{proof}
We have $V=\HOD$ and a bedrock in the constructible universe $L$. And the main result of \cite{HamkinsReitzWoodin2008:TheGroundAxiomAndVequalsHOD} provides a model of the ground axiom in which $V\neq\HOD$, which therefore exhibits the other prong of the dichotomy in a bedrock model.

Meanwhile, the bottomless model of Reitz is achieved by a set-support $\Ord$-product of adding a Cohen subset to every regular cardinal, which is homogeneous proper class forcing, so $\HOD^{V[G]}\of V$ and consequently it cannot be a ground of $V[G]$. So this model is in the $\necessary(V\neq\HOD)$ prong of the dichotomy, with no bedrock. 

For the remaining case, we want to provide a model of $V=\HOD$ which is bottomless. We start in a model $W$ of the continuum coding axiom \CCA, and perform a proper class $\Ord$ product of nonoverlapping self-encoding forcing $\P=\prod_\alpha\Q_\alpha$, where each $\Q_\alpha$ adds a Cohen subset of a regular cardinal, and then performs self-encoding of this (and itself) on a subsequent interval of cardinals. We may assume that the \CCA\ coding of $W$ takes place on an absolutely definable sequence of cardinals $\delta_\alpha$ that does not interfere with any segment of these various self-encoding forcing factors. Consider the resulting forcing extension $W[G]$, where $G$ is $W$-generic. We may view $G$ as the product of $g_\alpha\of\Q_\alpha$. Since the coding does not overlap or interfere, $W$ remains definable via the \CCA\ coding, and each $g_\alpha$ becomes ordinal-definable in $W[G]$, uniformly so, and so $G$ itself is definable in $W[G]$. It follows that $W[G]$ satisfies $V=\HOD$, but since this is a proper class product forcing, the model will be bottomless, since we can peel off the factors $\Q_\alpha$ one at a time.
\end{proof}

Let us next show that the contingent \HOD\ dichotomy does settle some features of set-theoretic geology. As we mentioned earlier, we know by \cite[theorem~68,70]{FuchsHamkinsReitz2015:Set-theoreticGeology}, using the main result of \cite{HamkinsReitzWoodin2008:TheGroundAxiomAndVequalsHOD}, that it is consistent with \ZFC\ that the mantle is distinct from the generic \HOD, but we should like to prove next that this requires the second prong of the contingent \HOD\ dichotomy.

\begin{theorem}
 In the $\possible(V=\HOD)$ prong of the contingent \HOD\ dichotomy, the mantle is equal to the generic \HOD.
  $$\Mantle=\gHOD$$
\end{theorem}

\begin{proof}
We know from \cite{FuchsHamkinsReitz2015:Set-theoreticGeology} that $\gHOD\of\Mantle$ in any model of \ZFC. But conversely, in the first prong of the contingent \HOD\ dichotomy, we know that $\HOD^{V[g]}$ is a ground model in every forcing extension $V[g]$. We may absorb the forcing into collapse forcing $V[g]\of V[h]$, which is homogeneous over $V$, leading to
 $$\HOD^{V[h]}\of\HOD^{V[g]}\quad\text{ and }\quad\HOD^{V[h]}\of\HOD\of V\of V[g].$$
Consequently, by the intermediate-model theorem for forcing (see \cite[corollary~15.43]{Jech:SetTheory3rdEdition}, \cite[Fact~11]{FuchsHamkinsReitz2015:Set-theoreticGeology}), it follows that $\HOD^{V[h]}$ is a ground of $V$. So the $\gHOD$ is intersecting models that are amongst the grounds of $V$, and so $\Mantle \of \gHOD$. So they are equal.

For a second, alternative argument, consider any forcing extension $V[g]$, and observe that because we are in the $\possible(V=\HOD)$ case of the dichotomy, we have both $V$ and $\HOD^{V[g]}$ being grounds of $V[g]$. By the downward directed grounds hypothesis applied in $V[g]$, therefore, there is a deeper ground $U\of V\intersect\HOD^{V[g]}$. In particular, $U$ is a ground of $V$ and also $U\of\HOD^{V[g]}$. So the mantle is contained in every such \HOD, and consequently $\Mantle\of\gHOD$, as desired.
\end{proof}

We should like also to observe that this property $\Mantle=\gHOD$, however, is not the same as the  $\possible(V=\HOD)$ prong of the contingent \HOD\ dichotomy. To see this, consider the model $V[G]$ used in the proof of \cite[theorem~66]{FuchsHamkinsReitz2015:Set-theoreticGeology}. This is the set-support $\Ord$-product of forcing that forces the \GCH\ at regular cardinal $\delta_\alpha$ or its negation, with the lottery sum of these two options, for an absolutely definable sequence of regular cardinals $\delta_\alpha$. This model has $\Mantle^{V[G]}=\gHOD^{V[G]}=\HOD^{V[G]}=V\ofneq V[G]$, so this is a case where $\Mantle=\gHOD$, and yet we are in the $\necessary(V\neq\HOD)$ prong of the contingent \HOD\ dichotomy.

Indeed, we might observe more generally that proper class homogeneous forcing will always place us into the second prong of the contingent \HOD\ dichotomy, since $\HOD^{V[G]}\of V$, and so $V[G]$ cannot have its \HOD\ as a ground. 

\section{The contingent \HOD\ dichotomy and iterated \HOD s}

We find it extremely natural to consider the iterated $\HOD^\alpha$ classes in the light of the contingent \HOD\ dichotomy. In the $\possible(V=\HOD)$ prong of the dichotomy, after all, the \HOD\ of any ground of $V$ is also a ground, and so in particular, $\HOD$ is a ground, and so $\HOD^\HOD$ is a ground, and also $\HOD^{\HOD^\HOD}$, and so forth. We define recursively that $\HOD^0=V$ and $\HOD^{n+1}=\HOD^{\HOD^n}$\!. In the first prong of the dichotomy, these are all grounds of $V$. 

One is naturally inclined to define $\HOD^\omega=\Intersect_n\HOD^n$, and indeed more generally $\HOD^{\alpha+1}=\HOD^{\HOD^\alpha}$ and $\HOD^\lambda=\Intersect_{\alpha<\lambda}\HOD^\alpha$ for limit ordinals. But there is a subtle issue that prevents this definition from being generally meaningful in \ZFC. Namely, we can often lack a uniform definition of $\HOD^\alpha$ and even of the $\HOD^n$ for finite iterates---the basic problem is that as we iterate the definition of \HOD, the resulting relativized definitions of $\HOD^n$ become more and more complicated, and it is not possible in \ZFC\ to provide a uniform definition of $\HOD^n$. So we cannot generally refer to all the classes $\HOD^\alpha$ in such a way so as to support the desired definition of the limit stage $\HOD^\lambda$. Harrington proved that this problem actually occurs and even $\HOD^\omega$ may not exist as a definable class---it might not even be \ZFC-amenable; see~\cite[section
7]{Zadrozny1983:IteratingOrdinalDefinability}. 
McAloon \cite{McAloon1971:ConsistencyResultsAboutOrdinalDefinability, McAloon1974:On-the-sequence-of-models-HOD_n} showed that even when it is a class, $\HOD^\omega$ may not be a model of \ZFC, but only $\ZF+\neg\AC$. 

Meanwhile, in some models of \ZFC, we do have a uniform treatment of the iterated \HOD s. See  \cite{HamkinsKirmayerPerlmutter2012:GeneralizationsOfKunenInconsistency} for the concept of having a \emph{uniform presentation} of the iterated $\HOD^\alpha$. Namely, the assertion ``$\HOD^\eta$ exists'' means that there is a parametrically definable sequence of classes $H^\alpha$ for $\alpha\leq\eta$, appearing as sections of a class $H\of(\eta+1)\times V$, such that $H^{\alpha+1}=\HOD^{H^\alpha}$ and $H^\lambda=\Intersect_{\alpha<\lambda}H^\alpha$ at limit stages. 

As we mentioned, the iterated \HOD\ sequence figures quite naturally in the first prong of the contingent \HOD\ dichotomy, since in this case \HOD\ is a ground, and therefore $\HOD^2$ as well and indeed every $\HOD^n$ is a ground. How far does it go? If $\HOD^\alpha$ exists and is a ground, then so is $\HOD^{\alpha+1}$. Perhaps one hopes to prove by induction that every $\HOD^\alpha$ is a ground.

Let us show at least that in the first prong of the contingent \HOD\ dichotomy, we avoid the Harrington-style obstacle regarding $\HOD^\omega$. 

\begin{theorem}
 Assume $\possible(V=\HOD)$, the first prong of the contingent \HOD\ dichotomy. Then $\HOD^\omega$ exists---we have a uniform presentation of $\HOD^\alpha$ for $\alpha\leq\omega$. Indeed, $\HOD^\eta$ exists as long as $\HOD^{\lambda+1}$ is a ground for all limit ordinals $\lambda<\eta$, if any. 
\end{theorem}

\begin{proof}
Assume $\possible(V=\HOD)$. It follows that $\HOD^W$ is a ground whenever $W$ is a ground, and so the iterated \HOD s continue to be grounds at successor stages. Suppose that $\HOD^{\alpha+1}$ is a ground for all $\alpha<\eta$ for some limit ordinal $\eta$. 

By the ground-model enumeration theorem \cite[theorem~12]{FuchsHamkinsReitz2015:Set-theoreticGeology}, there is a uniformly definable parameterized family $W_r$, for $r\in V$, such that every $W_r$ is a ground and every ground is $W_r$ for some $r$. Now the point is that 
we can express the relation $W_r=\HOD^n$ by asserting that there is a sequence of indices $\<r_0,\ldots,r_n>$ witnessing the \HOD\ iteration, that is, for which $W_{r_0}=V$ and $W_{r_{i+1}}=\HOD^{W_{r_i}}$ for each $i<n$. Furthermore, we can prove in $\ZFC+\possible(V=\HOD)$ that if $W_r=\HOD^n$, then there is a $r'$ such that $W_{r'}=\HOD^{n+1}$, since the \HOD\ of a ground is a ground. So we can definably pick out when an index $r$ is the index of a $\HOD^n$, and therefore $\HOD^\omega$ is definable as the intersection of these models. The main point is that the definitions of the iterated $\HOD^n$ do not become more and more complicated, as they do with a naive treatment in \ZFC, since instead of nesting the iterated definitions we are rather just picking out sequences of grounds by their indices and expressing instances of the same relation at each step, namely, that each next ground is the \HOD\ of the previous ground, which has bounded complexity in light of the ground-model enumeration theorem. So $\HOD^\omega$ exists. 

The same idea works more generally for $\HOD^\eta$, provided inductively that $\HOD^{\lambda+1}$ is a ground for all limit ordinals $\lambda<\eta$. The assertion $W_r=\HOD^\beta$ is expressed by the existence of a sequence of indices $\<r_\alpha\mid\alpha\leq\beta>$ for which $W_{r_0}=V$; $W_{r_{\alpha+1}}=\HOD^{W_{r_\alpha}}$;  $W_{\lambda_\alpha}=\Intersect_{\alpha<\lambda}W_{r_\alpha}$ for limit ordinals; and finally $r=r_\beta$. At bottom, the uniform presentation of grounds enables a uniform presentation of the $\HOD^\alpha$, as long as these remain grounds---and in the first prong of the contingent \HOD\ dichotomy, they will remain grounds at all successor stages.
\end{proof}

McAloon's model is obtained by forcing over $L$, so in his model every $\HOD^{\alpha+1}$ is an intermediate \ZFC\ model between the ground model $L$ and the universe $V=L[G]$, so these are all ground models, even though $\HOD^\omega$ is not a \ZFC\ model. Might we hope that this feature is true generally? 

\begin{question}
 If $\possible(V=\HOD)$, then is every $\HOD^{\alpha+1}$ a ground? Are the iterated $\HOD^\alpha$ coinitial in the grounds? Is the mantle the intersection of all $\HOD^\alpha$?
\end{question}

If we consider forcing in the context of \ZF\ models, it could be that $\HOD^\omega$ is a ground or indeed that every $\HOD^\alpha$ is a ground, even when they are not  \ZFC\ models. So it is natural to inquire under $\possible(V=\HOD)$ whether every $\HOD^\alpha$ is a \ZF-ground. 

We are unsure. We also wondered about another question:

\begin{question}
 Assume $\possible(V=\HOD)$ and there is a bedrock. Then, does every $W$ have its $\HOD^\alpha$ sequence eventually stabilizing? 
\end{question}

The idea would be that the bedrock would be a lower bound for the iterated $\HOD^\alpha$, and perhaps this causes them to stabilize. If they are all grounds, for example, then they would have to stabilize, since when there is a bedrock, there are only a set number of grounds, arising from Boolean subalgebras of the Boolean algebra realizing $V$ as a forcing extension of the bedrock.

We noticed that the following question also remains unsettled in our analysis, and we are unsure of the answer.

\begin{question}
 Is it consistent with \ZFC\ that $V\neq\HOD$ holds in all ground models, yet $V=\HOD$ holds in some forcing extension?
\end{question}

In light of theorem \ref{Theorem.Contingent-HOD-dichotomy-with-bedrock}, such a model would have to be bottomless, admitting no bedrock model, and every ground $W$ would sit atop a strictly descending chain of iterated $(\HOD^n)^W$, all being grounds but never stabilizing in a model of $V=\HOD$. But we don't know how to make this situation occur or to rule it out. 

\section{The contingent \HOD\ dichotomy and the maximality principles}

Let us turn now to the maximality principles, establishing that the contingent \HOD\ dichotomy is independent of the maximality principles \MP\ and $\MP(\R)$ over \ZFC.

\begin{theorem}\label{Theorem.V=HOD-and-MP-independence}
 The contingent \HOD\ dichotomy is independent of $\ZFC+\MP$. Indeed, the following theories are equiconsistent:
 \begin{enumerate}
  \item $\ZFC+\MP+(V=\HOD)$
  \item $\ZFC+\MP+\necessary(V\neq\HOD)$
  \item $\ZFC+\MP$
  \item $\ZFC$
  \item $\ZFC+(V_\delta\elesub V)$
  \item $\ZFC+(V_\delta\elesub V)+(V=\HOD)$
  \item $\ZFC+(V_\delta\elesub V)+\necessary(V\neq\HOD)$
 \end{enumerate}
\end{theorem}

The axiom $V_\delta\elesub V$ is the scheme expressed in the language of set theory augmented with a constant symbol for $\delta$ and the axioms: 
 $$\forall a\in V_\delta\bigl(\varphi(a)\iff\varphi^{V_\delta}(a)\bigr).$$
Such a cardinal $\delta$ is known as a \emph{correct} cardinal. We urge caution with this terminology, however, since the concept of being a correct cardinal is not expressible in the first-order language of set theory, except as a scheme, though it is expressible in Kelley-Morse set theory.

Another subtle point is that although $\ZFC+V_\delta\elesub V$ implies that every particular axiom of \ZFC\ holds in $V_\delta$, it does not imply $\Con(\ZFC)$ since no finitely many of those assertions are enough to get the whole theory. At the heart of this issue is the fact that $\Con(\ZFC)$ refers to the internal conception of \ZFC, which includes nonstandard axioms in any $\omega$-nonstandard model, whereas $V_\delta\elesub V$ concerns only the standard formulas of the metatheory.

\begin{proof}
The models arising in 1 and 2 establish the independence claim over \MP. The equiconsistency of theories 3, 4, and 5 is proved in~\cite{Hamkins2003:MaximalityPrinciple}. And clearly the consistency of any of the theories implies the consistency of theory 4, and theories 6 and 7 imply theory 5. There are several consistency implications remaining. 

$\bigl(\Con(5)\implies\Con(6)\bigr)$ 
If $V_\delta\elesub V$, then also $L_\delta\elesub L$, which also has $V=\HOD$, and so the $L$ of any model of theory 5 satisfies theory 6. One can also get theory 6 by the canonical (class) forcing of $V=\HOD$ over any model of theory 5, which will preserve the correctness of $\delta$, since if $V[G]_\delta\satisfies\varphi(a)$, then there is a name $\dot a\in V_\delta$ for $a$ and a condition $p\in G_\delta$ forcing $\varphi(\dot a)$ in $V_\delta$, which remains true in $V$, and so $\varphi(a)$ holds in $V[G]$. 

$\bigl(\Con(5)\implies\Con(7)\bigr)$ 
A similar argument works for theory 7, in light of the proof of theorem~\ref{Theorem.VneqHOD-can-be-necessary}. Start with $V_\delta\elesub V$ and then force with the class forcing to add a Cohen subset to every regular cardinal. The extension $V[G]$ satisfies $\necessary(V\neq\HOD)$ by theorem~\ref{Theorem.VneqHOD-can-be-necessary} and $V_\delta\elesub V$ is preserved just as in the previous paragraph. 

$\bigl(\Con(6)\implies\Con(1)\bigr)$ 
Assume $V_\delta\elesub V$ and $V=\HOD$. The main argument of~\cite{Hamkins2003:MaximalityPrinciple} shows that $\MP$ holds in the forcing extension $V[g]$ to collapse $\delta$ to become countable. And since we had $V=\HOD$ in $V$, it follows by theorem~\ref{Theorem.V=HOD-if-true-is-necessarily-forceable} that $V=\HOD$ is a switch, so there is a further forcing extension $V[g][h]$ in which $V=\HOD$ holds. Furthermore, it is a general observation about the maximality principle that whenever it holds, then it also holds in every forcing extension. So $V[g][h]$ is a model of $V=\HOD$ plus $\MP$, as desired. 

$\bigl(\Con(7)\implies\Con(2)\bigr)$ Assume $V_\delta\elesub V$ plus $\necessary(V\neq\HOD)$. After collapsing $\delta$ to become countable, we achieve $\MP$ in the forcing extension $V[g]$. And since this is a set-forcing extension, we still have $\necessary(V\neq\HOD)$, as desired. 
\end{proof}

A similar argument will show that the contingent \HOD\ dichotomy is independent of the boldface maximality principle $\MP(\R)$. 

\begin{theorem}\label{Theorem.V=HOD-and-MPtilde-independence}
 The contingent \HOD\ dichotomy is independent of $\ZFC+\MP(\R)$. More specifically, the following theories are equiconsistent:
 \begin{enumerate}
  \item $\ZFC+\MP(\R)+(V=\HOD)$.
  \item $\ZFC+\MP(\R)+\necessary(V\neq\HOD)$.
  \item $\ZFC+\MP(\R)$. 
  \item $\ZFC+\Ord$ is definably Mahlo.
  \item $\ZFC+(V_\delta\elesub V)+(\delta\text{ is inaccessible})$.
  \item $\ZFC+(V_\delta\elesub V)+(\delta\text{ is inaccessible})+(V=\HOD)$.
  \item $\ZFC+(V_\delta\elesub V)+(\delta\text{ is inaccessible})+\necessary(V\neq\HOD)$.
 \end{enumerate}
\end{theorem}

The axiom ``$\Ord$ is definably Mahlo'' is the scheme asserting that every definable club class of ordinals contains a regular cardinal, or equivalently, that every such class club contains an inaccessible cardinal. 

\begin{proof}
The models arising in 1 and 2 establish the independence claim over $\MP(\R)$. The theories in 3, 4, and 5 are proved equiconsistent in~\cite{Hamkins2003:MaximalityPrinciple}. Clearly the theories in 1 and 2 imply theory 3, and theories 6 and 7 imply theory 5. So just a few consistency implications remain. 

$\bigl(\Con(5)\implies\Con(6)\bigr)$  
Assume $\delta$ is inaccessible and $V_\delta\elesub V$. Then this is also true if we move to $L$, where we have $V=\HOD$, as desired. Alternatively, just as in theorem~\ref{Theorem.V=HOD-and-MP-independence}, we can also achieve $V=\HOD$ by the canonical class forcing of $V=\HOD$, which preserves $V_\delta\elesub V$ and the inaccessibility of $\delta$.   

$\bigl(\Con(5)\implies\Con(7)\bigr)$  
Assume $\delta$ is inaccessible and $V_\delta\elesub V$, and perform class forcing $V\of V[G]$ with the Easton-support product that adds a Cohen subset to every regular cardinal. This preserves $V_\delta\elesub V$ by the argument we gave earlier, and it also preserves the inaccessibility of $\delta$ by the usual Easton argument. And it forces $\necessary(V\neq\HOD)$ by theorem~\ref{Theorem.VneqHOD-can-be-necessary}. So we have achieved theory 7 in the forcing extension $V[G]$, as desired.

$\bigl(\Con(6)\implies\Con(1)\bigr)$  
Assume $V_\delta\elesub V$, where $\delta$ is inaccessible, plus $V=\HOD$. The main result of~\cite{Hamkins2003:MaximalityPrinciple} shows that if we force with the \Levy\ collapse of $\delta$, that is, with $\text{Coll}(\omega,{<}\delta)$, then we get a forcing extension $V[g]$ with the $\MP(\R)$. Since $V=\HOD$ holds in $V$, we know that it is forceable over $V[g]$, and furthermore, we can perform the self-encoding forcing entirely above $\delta$, which will not add reals. The result is a forcing extension $V[g][h]$ in which $V=\HOD$. Since we did not add any real numbers, we retain every instance of the $\MP(\R)$ scheme from $V[g]$, and so $V[g][h]$ also satisfies $\MP(\R)$. So theory 1 holds in $V[g][h]$, as desired.

$\bigl(\Con(7)\implies\Con(2)\bigr)$  
Assume $V_\delta\elesub V$, where $\delta$ is inaccessible, plus $\necessary (V\neq\HOD)$. After the \Levy\ collapse of $\delta$, we get $\MP(\R)$ in the extension $V[g]$. And since $\necessary(V\neq\HOD)$ held in $V$, we still have it in $V[g]$, as desired.
\end{proof}

For completeness we can also mention the dual result that the maximality principles are independent of $V=\HOD$ over \ZFC.  We have already provided models of $V=\HOD$ plus the boldface maximality principle $\MP(\R)$ and hence \MP\ as well. What remains is to find a model of $\ZFC+V=\HOD$ in which \MP\ fails. But this is quite easy, for this is the situation in the constructible universe $L$. The summary conclusion, therefore, is that $V=\HOD$ and the maximality principles are independent of each other over \ZFC.

Finally, regarding the necessary maximality principle $\necessary\MP(\R)$, which asserts that $\MP(\R)$ holds in every forcing extension, Wojciech Wo{\l}oszyn has recently announced on Twitter (take it as you will) a result that if the theory ``$\ZF+\AD^++\AD_\R+\Theta$ is regular'' is consistent, then so is $\ZFC+\necessary\MP(\R)+(V\neq\HOD)$.

\section{The contingent \HOD\ dichotomy compared with Woodin's \HOD\ dichotomy}\label{Section.Contingent-HOD-dichotomy-versus-HOD-dichotomy}

Let us finally discuss the precise relation of the contingent \HOD\ dichotomy of this paper and the well-known \HOD\ dichotomy of W. Hugh Woodin, as discussed in \cite{Woodin_Davis_Rodriguez_2012}. Although the two dichotomies both fall under the informal description that $V$ is either close to \HOD\ or far from \HOD, nevertheless the dichotomies are not identical and they do not partition set-theoretic possibility in exactly the same manner. Woodin proposes his dichotomy under the assumption that there is an extendible cardinal:

\begin{theorem}[Woodin's HOD dichotomy
  \cite{Woodin_Davis_Rodriguez_2012}]
\label{thm:woodin-hod-dichotomy}
Suppose that $\delta$ is an extendible cardinal. Then exactly one of the
following alternatives holds.
\begin{enumerate}
\item Every singular cardinal $\lambda>\delta$ is singular in
$\mathrm{HOD}$ and $(\lambda^+)^{\mathrm{HOD}}=\lambda^+$.
\item Every regular cardinal $\kappa>\delta$ is $\omega$-strongly
measurable in $\mathrm{HOD}$. 
\end{enumerate}
\end{theorem}

Goldberg \cite{Goldberg2024StronglyCompact} subsequently shows that the large-cardinal hypothesis can be reduced to a strongly compact cardinal. The \emph{\HOD\ Hypothesis} is the assertion that there is a proper class of regular cardinals that are not $\omega$-strongly measurable in $\HOD$. The \emph{\HOD\ Conjecture} asserts that ZFC + ``there is an extendible cardinal'' proves the HOD Hypothesis \cite{Woodin_Davis_Rodriguez_2012}.

We should like to point out the first prong of our contingent \HOD\ dichotomy is contained within the first prong of Woodin's \HOD\ dichotomy. This is simply because if \HOD\ is a ground, as in the first prong of our contingent \HOD\ dichotomy, then it is correct about all regular cardinals above the size of the relevant forcing; therefore \HOD\ will be correct about all sufficiently large successor cardinals, which will consequently not be measurable in \HOD, and this prevents the second prong of Woodin's dichotomy, so we must be in his first prong. In short, our first prong is a refinement of Woodin's first prong.

But the two dichotomies are not identical. To see this, consider the central model $V[G]$ appearing in the main theorem of \cite{HamkinsReitzWoodin2008:TheGroundAxiomAndVequalsHOD}, which has $V\neq\HOD$ plus the ground axiom. By theorem \ref{Theorem.Contingent-HOD-dichotomy-with-bedrock}, this places it into the $\necessary(V\neq\HOD)$ prong of the contingent \HOD\ dichotomy of this paper. And yet, that model proceeds by an $\Ord$-iteration that forces to add a Cohen subset at each regular cardinal (or at an absolutely definable class of such cardinals), and this is forcing that will preserve every extendible cardinal by the usual master-condition lifting arguments (for example, \cite[theorem~2.2]{Tsaprounis2013:OnExtendibleCardinalsAndGCH}). 
This forcing is mild and will preserve the first prong of Woodin's dichotomy. In summary, this model appears in the first prong of Woodin's \HOD\ dichotomy, but in the second prong of our contingent \HOD\ dichotomy, showing that the dichotomies are not the same.

In a sense, since Woodin has conjectured that the second prong of his dichotomy does not occur, we find it helpful to offer a natural refinement of the first prong of his dichotomy, which nevertheless still carves the set-theoretic multiverse at its joints, explicating how the set-theoretic universe might interact with its \HOD.

We have already mentioned in theorem \ref{Theorem.Large-cardinals-consistent-with-both-prongs} that large cardinals are generally relatively consistent with both prongs of the contingent \HOD\ dichotomy. Meanwhile, theorem \ref{Theorem.Goldberg-ultrapower-axiom} shows that under the ultrapower axiom with a supercompact, \HOD\ is a ground, placing us in the first prong of the contingent \HOD\ dichotomy. 

In contrast, Aguilera, Bagaria, and L\"ucke recently introduced \emph{exacting} and \emph{ultraexacting} cardinals \cite[Definitions~2.4 and~3.3]{AguileraBagariaLuecke2025}, making a connection with Woodin's \HOD\ dichotomy. Exacting cardinals imply $V\neq\mathrm{HOD}$, and ultraexacting cardinals form a strengthening of exacting cardinals. Let us draw out some relevant observations concerning how their results interact with our contingent \HOD\ dichotomy.

\begin{observation}\ 
\begin{enumerate}
\item If there is an extendible cardinal below an exacting cardinal,
then the second prong of the contingent HOD dichotomy holds.
\item Relative to the consistency of an $\mathrm{I}2$ embedding, it is
consistent that there is an exacting cardinal below an extendible
cardinal and  the first prong of the contingent HOD dichotomy holds.
\item Relative to the consistency of $\mathrm{I}0$, it is consistent
that there is an ultraexacting cardinal and the first prong of the contingent HOD dichotomy holds.
\end{enumerate}
\end{observation}

\begin{proof}
For (1), Aguilera--Bagaria--L\"ucke prove that an exacting cardinal
above an extendible cardinal implies that the HOD Hypothesis
fails \cite[Conclusion~6.8]{AguileraBagariaLuecke2025}. So by the observation above the second prong
of the contingent HOD dichotomy follows. 

For (2), Aguilera--Bagaria--Goldberg--L\"ucke show, relative to an
$\mathrm{I}2$ embedding, the consistency of the HOD Hypothesis together
with an exacting cardinal below an extendible cardinal
\cite[Corollary~5.6]{AguileraBagariaGoldbergLuecke2025}. Their proof
starts from a transitive model of
\[
 \text{ZFC}+V=\mathrm{HOD}+\text{Vop\v{e}nka's Principle}
\]
and obtains the desired configuration by set forcing. Hence the final
model has a ground satisfying $V=\mathrm{HOD}$, and so it lies in the
first prong.

For (3), in \cite[Theorem~4.5]{AguileraBagariaGoldbergLuecke2025} the model has the form $ N=L(V_{\lambda+1})[G]$  satisfying ZFC and ``$\lambda$ is ultraexacting''. In $N$, we can code $\langle V_{\lambda+1},G\rangle$ by a set of ordinals $b$. Then
$N=L[b]$, and consequently $N$ satisfies the first prong by theorem \ref{Theorem.L[A]-is-in-first-prong}.
\end{proof}

The following seems to be open.

\begin{question}
\label{question:extendible-above-ultraexacting-second-prong}
Does the existence of an extendible cardinal above an ultraexacting
cardinal imply the second prong of the contingent HOD dichotomy?
\end{question}

Aguilera--Bagaria--Goldberg--L\"ucke ask whether this configuration
disproves the HOD Hypothesis
\cite[Question~7.1]{AguileraBagariaGoldbergLuecke2025}. A positive
answer to their question would give a positive answer to
Question~\ref{question:extendible-above-ultraexacting-second-prong}. 

\printbibliography

\end{document}